\documentclass[11pt,twoside]{article}
\usepackage{graphicx} 
\usepackage{ulem}
\usepackage{amssymb}
\usepackage[mathscr]{eucal}
\usepackage{eufrak}
\usepackage{amsmath,amsthm}
\usepackage{physics}
\usepackage{mathrsfs}
\usepackage{lgreek}
\usepackage[shortlabels]{enumitem}
\usepackage[colorlinks=true,
   urlcolor=blue,           filecolor=green,      
   citecolor=green,      
   linkcolor=red,           bookmarks=true,
  unicode,
   plainpages=false,   ]{hyperref}
\usepackage{color}
\usepackage{comment}
\font\bg=cmbx10 scaled 1300

\definecolor{royalblue}{rgb}{0,0,0.128}
\def\bp{\begin{proof}}
\def\ep{\end{proof}}
\def\n{\nabla}

\def\dm{|\hskip-0.05cm|}

\def\OO{\Omega}
\def\displ{\displaystyle}

\def\VS{\vspace{6pt}\\\displ }

\def\R{\Bbb R}
\def\N{\Bbb N}
\def\à{à}

\def\vep{\varepsilon}

\def\be{\begin{equation}}
\def\ba{\begin{array}}
\def\ea{\end{array}}
\def\ee{\end{equation}}
\def\vs1{\vspace{1ex}}

\def\Ã©{\'{e}}
\def\Ãš{\`{e}}

\newtheorem{lemma}
{\bf Lemma} 
\font\sc=cmcsc10
\newtheorem{defi}
{\bf Definition} 
\newtheorem{tho}
{\bf Theorem} 
\newtheorem{rem}
{\sc Remark} 
\newtheorem{coro}
{\bf Corollary}

\title{\bg An $L^p$-theory for global weak solutions to the Navier-Stokes equations in exterior domains}
\author{\sc  Filippo Palma
\thanks{Dipartimento di Matematica e Fisica,  
Universit\`{a} degli 
Studi della Campania
``L. Vanvitelli'', via Vivaldi 43, 81100 \null\hskip0.55cmCaserta,
 Italy.\newline\null\hskip0.55cm
 filippo.palma@unicampania.it
\newline\null\hskip0.55cm The  research activity  is performed under the
auspices of   GNFM-INdAM. }}
\begin{document}
\markboth{\footnotesize\rm    F. Palma} {\footnotesize\rm An $L^p$-theory for global weak solutions to the Navier-Stokes equations in exterior domains
}
\maketitle
\begin{abstract}
    In this paper, we consider the Navier-Stokes initial boundary value problem in exterior domains with initial data in the Lebesgue space $L^p$, $p\in (2,3)$, and show the existence of a global weak solution. For the quoted solution, we also furnish a structure theorem. In particular, we see that the weak solution becomes regular after a certain instant of time and it is also regular a. e. in time. Although a general $L^p$-theory for local strong/mild solutions is well-established in the literature, a corresponding theory for global weak solutions in exterior domains seems to be new. Of course, our results hold in the particular cases of the Cauchy problem and the initial boundary value problem in a half-space.
\end{abstract}
\noindent
\section{Introduction}
Let $\Omega\subseteq\R^3$ denote an exterior domain, whose boundary for simplicity is supposed to be sufficiently smooth, the whole space, or a half-space. We consider the Navier-Stokes initial boundary value problem (initial value problem if $\Omega=\R^3$):
\begin{equation}\label{NS-IBVP}
\begin{cases}
    u_t-\Delta u+u\cdot\n u+\n\pi_u=0\,,\quad\,\,\,\,\,\, \text{in }(0,T)\times\Omega\,,\\
   \hskip3.27cm \n\cdot u=0\,,\quad\,\,\,\,\,\,\text{in }(0,T)\times\Omega\,,\\
   \hskip3.05cmu(t,x)=0\,,\quad \,\,\,\,\,\,\text{in }(0,T)\times\partial\Omega\,,\\
   \hskip2.3cm (u(0,x),\varphi)=(u_0(x),\varphi)\,,\,\text{for all }\varphi\in \mathscr C_0(\Omega)\,.
\end{cases}
\end{equation}
The initial condition \eqref{NS-IBVP}$_4$ is considered in a weak formulation since we are going to investigate weak solutions corresponding to initial data that do not need to be necessarily weakly divergence-free. This kind of initial condition was already employed in \cite{Mar-L1} in the context of the Stokes problem. However, in the case of a weakly divergence-free initial datum, the problem can be reformulated in the usual way. \par The aim of this work is to show the existence and investigate the regularity of weak solutions to problem \eqref{NS-IBVP} corresponding to initial data in $L^p(\OO)$, $p\in (2,3)$. In the particular case of the Cauchy problem, a first study in this direction was carried out by Calder\'on in \cite{Cald}. Employing a decomposition of the initial datum as the sum of a ``small" $L^3$ part and a ``big" $L^2$ part, he reduced the Cauchy problem for the Navier-Stokes equations to two nonlinear problems investigated by different techniques. However, he obtained a solution satisfying a very weak formulation of the Navier-Stokes problem and no results related to the partial regularity of the quoted solution are furnished. Our goal is, on the one hand, to fill this gap, by obtaining a local integrability condition for the first-order spatial derivatives as long as proving a structure theorem for our weak solution. On the other hand, we employ more flexible and immediate arguments that allow to get our results in more general domains, such as exterior domains and the half-space. \par In a different context, a decomposition argument was employed by Maremonti \cite{Mar-Ferrara, Mar-L3Lor} and Maremonti and Shimizu \cite{Mar-Shi-2D, Mar-Shi-3D}. However, in these last works, the decomposition employed is related to the equations of the fluid motion that are reduced to the study of different subproblems, among which only one has a nonlinear character. The initial datum is not decomposed and is considered as the initial datum for a suitable Stokes problem. \par On the other hand, employing a Calder\'{o}n-like decomposition, in \cite{Pop} Popkin obtained, in the case of the Cauchy problem, global weak solutions for data in supercritical Besov spaces. \par Since the pioneering work of Leray \cite{Ler}, it is well-known that for all data in $L^2$ there exists at least one global weak solution to problem \eqref{NS-IBVP}. Moreover, the weak solution enjoys a partial regularity property known as \textit{Th\'eor\`eme de Structure}. Leray dealt with the Cauchy problem and his results were extended to the more general case of the initial boundary value problem by several authors. \par Leray's choice of the initial datum is performed in order to get a solution having finite energy for all times. If we do not require an initial datum that belongs to some subspace of $L^2(\OO)$, the finiteness of the kinetic energy of the fluid cannot be guaranteed. However, the study of the Navier-Stokes IBVP for initial data that do not have finite energy led to many interesting results. \par One of the first contributions in this direction is the well-known paper of Fabes, Jones and Rivi\`ere \cite{FJR}. In their work, the authors prove the existence and uniqueness of a local solution to the integral formulation of the Navier-Stokes Cauchy problem for any initial datum in $L^q(\R^3)$, $q>3$. Moreover, the quoted solution satisfies a very weak formulation of the Navier-Stokes Cauchy problem. Later, the result was extended to more general domains by Fabes, Lewis and Rivi\`ere \cite{FLR}. As in the previous work, the main idea of the proof consists in considering the integral formulation of the Navier-Stokes problem and employing the theory of hydrodynamic potentials. \par Another relevant tool that leads to a local mild solution to the Navier-Stokes problem for $L^q$ initial data is the analyticity of the Stokes semigroup on $L^q$ Lebesgue spaces \cite{Giga86, Giga-Miy}. The approach proposed by Giga in \cite{Giga86} works as well in more general domains and is extendable to other evolution problems, while the solution obtained by Miyakawa and Giga in \cite{Giga-Miy} enjoys further regularity properties compared to that obtained by Fabes, Lewis and Rivi\`ere and is called a \textit{strong solution}. \par {However, in all the previously quoted works, if $q>3$ the solution can be guaranteed only locally in time for general initial data in $L^q$.}  \par In the special case of an initial datum belonging to the scaling invariant space $L^3(\R^3)$, and satisfying a suitable smallness condition, Kato \cite{Kato} showed the existence of a unique strong solution that is global in time, see also \cite{Giga86, Iwa, Mar-DCDS} for some extensions of Kato's result to more general domains. \par A first result on the existence of global weak solutions to the Navier-Stokes initial boundary value problem corresponding to an initial datum in the scaling invariant Lorentz space $L^{3,\infty}(\OO)$ was obtained by Maremonti in \cite{Mar-L3Lor}, see also \cite{BSS}. The result holds for data of arbitrary size and a structure theorem for the weak solution is also furnished. A result for data in $L^3(\R^3)$ was obtained by Seregin and Sverak in \cite{Ser-Sve}, but no partial regularity was proved for the weak solution. In the quoted papers, a decomposition argument is also employed, but the decomposition does not involve the initial datum, which is assumed as the initial datum for a suitable Stokes problem. An analogous decomposition was also employed, in the context of the Cauchy problem, in \cite{Brad-Hud}, where a global weak solution is constructed for data in $L^{p,\infty}(\R^3)$, $p\in (2,3)$. \par For initial data that do not have necessarily a finite-energy and do not belong to scaling-invariant spaces, a result for local regular solutions is due to Heywood \cite{Hey}. In particular, he assumes an initial datum having a square-summable gradient. Under the same condition on the initial datum, Maremonti obtained a global weak solution, as long as a special structure theorem \cite{Mar-Ferrara}. However, the aforementioned structure theorem does not allow to deduce a time interval of regularity of the kind $(\theta,\infty)$, $\theta\geq 0$. This fact does not allow for the investigation of the long-time behavior of the quoted weak solution, as uniqueness cannot be guaranteed even for large times. This aspect deserves accurate further investigation.\par In the context of the 3D initial boundary value problem for nondecaying initial data, a global weak solution was obtained by Maremonti and Shimizu in \cite{Mar-Shi-3D}, where the initial datum is considered in $L^\infty(\OO)\cap \widehat W^{1,q}(\OO)$, $q>3$. \par In the special case of the Cauchy problem, many results for weak solutions corresponding to rough initial data were obtained in the last years. We mention, among others, the work by Lemari\'e-Rieusset on the so called ``Local-Energy Solutions" \cite{Lem-Rie} and the work by Kwon and Tsai on global weak solutions for data with ``slowly-decaying" oscillations \cite{Kwon-Tsai}. \par For further interesting results we also refer to \cite{Tsai-intermediate, Brad-Kuc, Brad-Kuc-Oza,  Brad-Kuc-Tsai}. \par We briefly introduce the notation employed throughout the paper.
We set
\[
\mathscr C_0(\OO):=\{\varphi\in C_0^\infty(\OO)\,:\,\n\cdot\varphi=0\}\,.
\]
By the symbol $\mu$ we refer to the Lebesgue measure on $\R^3$. We denote by $L^p(\OO)$, $W^{m,p}(\OO)$ and $\widehat W^{m,p}(\OO)$, $m\in\N$, $p\in \R^+$, respectively, the Lebesgue, Sobolev and homogeneous Sobolev spaces.\par 
We define the spaces $J^p(\OO)$ and $J^{m,p}(\OO)$ as the closures of $\mathscr C_0(\OO)$ with respect to the metrics $\dm \cdot\dm_p$ and $\dm \cdot\dm_{m,p}$, respectively. It is known that, for all $p\in (1,\infty)$, the following decomposition holds:
\[
L^p(\OO)=J^p(\OO)\oplus G^p(\OO)\,,
\]
where $G^p(\OO)$ denotes the subspace of $L^p(\OO)$ of functions of the kind $f=\n h$, $h\in W^{1,p}_{loc}(\OO)$ with $\n h\in L^p(\OO)$. Moreover, we are able to define a continuous operator
\[
\mathbb P_p:L^p(\OO)\to J^p(\OO)\,,
\]
known as Helmholtz projection. We also recall that there exists $c=c(p)>0$ such that for all $H\in L^p(\OO)$, $H=U\oplus \n K$, with $U\in J^p(\OO)$, we have
\be\label{HD-ineq}
\dm U\dm_p+\dm \n K\dm_p\leq c(p)\dm H\dm_p\,.
\ee
For more details see \cite{Ga:Book, Sol}.\par 
By the symbol $L^{p,q}(\OO)$, $p,q\in [1,\infty]$, we refer to the Lorentz spaces and we denote by $\dm\cdot\dm_{pq}$ the corresponding norms. Given $-\infty\leq a < b\leq+\infty$, we denote by  $C(a,b;X)$ and $L^p(a,b;X)$, $X$ Banach space, the Bochner spaces. We will usually employ the concise notation $(u,\phi)$ to denote the integral 
\[
\int_\OO u\cdot\phi\,dx\,.
\]
We refer to the trace space of $W^{m,p}(\OO)$ with the notation $W^{m-\frac1p,p}(\partial\OO)$ If $p=2$, we might employ the alternative notation $H^{m-\frac12}(\partial\OO)$. For a complete overview on function spaces we refer the reader to \cite{Adams}, while for the spaces of the hydrodynamics we quote \cite{Ga:Book}. \par Employing the notation in \cite{Ga:Book}, we denote by $B_R$ a ball centered in the origin with radius $R$, and we set $\OO_R=\OO\cap B_R$ and $\OO^R=\OO\setminus B_R$. \par In the following, by $c$ we mean a generic constant, whose value may vary from time to time and is inessential for our purposes. \par We give the following notion of global weak solution to problem \eqref{NS-IBVP}. We denote by \newline $\displ \gamma: W^{1,2}(\OO)\to H^{\frac12}(\partial\OO)$ the Gagliardo trace operator.
\begin{defi} \label{def: WS-NS}
    {\sl Let $u_0\in L^p(\OO)$, $p\in (2,3]$. A field $u:(0,\infty)\times\OO\to \R^3$ is said a weak solution to problem \eqref{NS-IBVP} if, for all $T>0$,
    \begin{description}
        \item[a)] $u\in L^{2,\infty}(0,T;W^{1,2}_{loc}(\OO))\cap L^4(0,T;L^3(\OO))$;
        \item[b)] $$\ba{c}\displ \int_0^T \big[(u,\phi_t)-(\n u,\n \phi)+(u\cdot\n \phi,u)\big]\,d\tau+(u_0,\phi(0))=0\,,\VS \text{for all }\phi\in C_0^\infty([0,T)\times\OO)\,,\,\,\text{with }\n\cdot\phi=0\,\,\text{for all }t\in [0,T)\,;\ea$$
        \item[c)] $(u(t),\n\psi)=0$, for all $\psi\in C_0^\infty(\OO)$ and $t\in (0,T)$;
        \item[d)] $\gamma (u(t))=0$ in $H^\frac12(\partial\OO)$ for a. a. $t\in (0,T)$ and $\displ\lim_{t\to 0}(u(t),\varphi)=(u_0,\varphi)$, for all $\varphi\in \mathscr C_0(\OO)$.
    \end{description}}
\end{defi}
We also give the notion of a regular solution.
\begin{defi}\label{def: SOLR-NS}
    {\sl Let $u_0\in L^3(\OO)$ and $u$ be a weak solution to problem \eqref{NS-IBVP} in the sense of Definition\,\ref{def: WS-NS}. $u$ is said a regular solution on $(0,T)\times\OO$ if
    \[
    \ba{c}
    \displ\eta >0\,,\,\, u\in C([0,T);J^3(\OO))\cap L^\infty(\eta,T;J^{1,3}(\OO)\cap W^{2,3}(\OO))\,,\, u_t\in L^\infty(\eta,T;L^3(\OO))\,,
    \ea
    \]
    and there exists a field $\pi_u:(0,T)\times\OO\to \R$ such that
    \[
    \displ\eta >0\,,\,\,\n\pi_u\in L^\infty(\eta,T;L^3(\OO))\,.
    \]}
\end{defi}
\begin{rem}
    {\rm
    \begin{description} 
    \item[i.] Actually, if the initial datum $u_0$ is not weakly divergence-free, the continuity in $t=0$ stated in Definition\,\ref{def: SOLR-NS} has to be meant in the following way: \[\lim_{t\to 0}\dm u(t)-\mathbb P_3u_0\dm_3=0\,;\]
    \item[ii.] The notions of weak and regular solution are analogous to those given in the context of the $L^{3,\infty}$-theory for global weak solutions due to Maremonti \cite{Mar-L3Lor}. Actually, the surprising aspect of our results is that a global weak solution corresponding to an $L^p$ initial datum, $p\in (2,3)$ behaves, a. e. in time and after a certain instant of time, as an $L^3$ regular solution. Moreover, we are able to show the existence of a global weak solution having the same regularity properties of a weak solution corresponding to an $L^3$ initial datum. 
    \end{description}}
\end{rem}
Our aim is to show that for all data in $L^p(\OO)$, $p\in (2,3)$ there exists a global weak solution to problem \eqref{NS-IBVP} in the sense of Definition\,\ref{def: WS-NS}. Moreover, for the quoted weak solution, we prove partial regularity in the spirit of Leray's \textit{Th\'eor\`eme de Structure}. \par In order to reach our goal, we partially follow an idea employed by Calder\'on in \cite{Cald}. We employ a decomposition of the initial datum as the sum $u_0=v_0+w_0$, with $v_0\in L^3(\OO)\cap L^p(\OO)$ and $w_0\in L^2(\OO)\cap L^p(\OO)$, which is possible by virtue of Lemma\,\ref{Dec-Lp}. Hence, the difference between our decomposition and the one proposed by Calder\'on consists in the fact that we do not require a weak divergence-free condition, which is anyway always possible, to be preserved under the decomposition. Our decomposition seems to be more immediate and it holds in general domains. In the special case of the Cauchy problem, Calder\'on proved the existence of a solution to the integral formulation of the Navier-Stokes problem satisfying a very weak formulation of problem \eqref{NS-IBVP}, \cite{Cald}. In our case, we get a local integrability property for the first order spatial derivatives, an existence result in more general domains, and last but not least, a structure theorem for our weak solution. \par We construct our weak solution $u$ as the sum $u=v+w$, where $v$ is a global regular solution to 
\begin{equation}\label{NS-aux-intro}
    \begin{cases}
    v_t-\Delta v+v\cdot\n v+\n\pi_v=0\,,\quad\,\,\,\,\,\,\,\,\,\,\,\quad \text{in }(0,T)\times\Omega\,,\\
   \hskip3.16cm \n\cdot v=0\,,\quad\,\,\,\,\,\,\quad\,\,\,\,\,\text{in }(0,T)\times\Omega\,,\\
   \hskip2.91cm v(t,x)=0\,,\quad \,\,\,\,\,\,\quad\quad\text{in }(0,T)\times\partial\Omega\,,\\
   \hskip2.14cm (v(0,x),\varphi)=(v_0(x),\varphi)\,,\,\text{for all }\varphi\in\mathscr C_0(\Omega)\,,
\end{cases}
\end{equation}
 while $w$ is a global weak solution to 
 \begin{equation}\label{NS-pert-intro}
      \begin{cases}
        w_t-\Delta w+w\cdot\n w+w\cdot\n v+v\cdot\n w+\n\pi_{w}=0\,,\,\, \text{in }(0,T)\times\Omega\,,\\
   \hskip3.16cm \n\cdot w=0\,,\quad\,\,\,\,\,\,\quad\,\,\,\,\,\,\quad\qquad\qquad\,\,\,\,\text{in }(0,T)\times\Omega\,,\\
   \hskip2.89cm w(t,x)=0\,,\quad\quad \,\,\,\,\,\,\quad\quad\qquad\qquad\,\,\,\,\,\text{in }(0,T)\times\partial\Omega\,,\\
   \hskip2.14cm (w(0,x),\varphi)=(w_0,\varphi)\,,\,\,\text{for all }\varphi\in\mathscr C_0(\Omega)\,.
    \end{cases}
 \end{equation}
 The well-posedness results for problems \eqref{NS-aux-intro}-\eqref{NS-pert-intro} are established, respectively, in \autoref{sec: auxprob} and \autoref{sec: pertprob}. \par
Our main result is the following.
\begin{tho}\label{thm: MR}
    {\sl Let $p\in (2,3)$ and $u_0\in L^p(\OO)$. Then there exists a weak solution $u:=v+w$ to problem \eqref{NS-IBVP} in the sense of Definition\,\ref{def: WS-NS}, where $(v,\pi_v)$ is a regular solution to problem \eqref{NS-aux-intro} and $w$ is a weak solution to problem \eqref{NS-pert-intro}. Moreover, there exist a sequence of open intervals of time $\{(\theta_l,T_l)\}_{l\in\N}$ and $\theta\geq 0$ such that 
    \[
\mu\bigg((0,\infty)\setminus\big(\underset{l\in\N}\cup (\theta_l,T_l)\cup [\theta,\infty)\big)\bigg)=0
    \]
    and $u$ is regular, in the sense of Definition\,\ref{def: SOLR-NS}, on $[\theta,\infty)$ as long as on $(\theta_l,T_l)$, for all $l\in\N$.}
\end{tho}
Since we are able to furnish a weak solution to problem \eqref{NS-IBVP}, which is smooth for $t\geq \theta$, then, in accord with the results by Kato \cite[Thereom 4]{Kato} and the semigroup properties of the Stokes operator, we can establish the following result on the asymptotic properties in time for our weak solution $u$.
\begin{coro}\label{cor: LTB-sol}
    {\sl Let $u$ be the solution stated in Theorem\,\ref{thm: MR}. Then
    \[
    \ba{c}
    \dm u(t)\dm_q\leq c \dm u(\theta)\dm_pt^{-\frac32(\frac1p-\frac1q)}\,,\,\,\text{for }t>\theta\,,\,\,\text{for all }q\in[p,\infty]\,,\VS\begin{cases} \dm \n u(t)\dm_q\leq c \dm u(\theta)\dm_pt^{-\frac32(\frac1p-\frac1q)-\frac12}\,, \,\,\text{for }t>\theta,\,\text{for all }q\in[p,\infty),\,\text{if }\OO=\R^3\,\text{or }\OO=\R^3_+\,, \\ \dm \n u(t)\dm_q\leq c \dm u(\theta)\dm_pt^{-\frac32(\frac1p-\frac1q)-\frac12}\,,\,\,\text{for }t>\theta,\,\text{for all }q\in [p,3],\,\text{if }\OO\,\,\text{is an exterior domain}. \end{cases}
    \ea
    \]}
\end{coro}
\begin{rem}
    {\rm Concerning Theorem\,\ref{thm: MR}, we want to stress that we are able to furnish an upper bound to the starting instant of time $\theta$ for the global regularity of our weak solution (see the proof of Theorem\,\ref{Struct-w}). An explicit upper bound for the starting instant of regularity is also known for the Leray weak solution. However, this bound was not given in the context of the weak solution corresponding to an initial datum in $L^{3,\infty}(\OO)$, see \cite{Mar-L3Lor}. We claim that our arguments allow us to get a similar explicit bound even in the $L^3$ or $L^{3,\infty}$ theory for global weak solutions. This aspect will be discussed in a forthcoming note. \par A future step consists in the investigation of global weak solutions for initial data in weighted Lebesgue spaces. For these special classes of initial data, ``small" perturbations of steady-state solutions in the whole space were investigated in \cite{Gal-Mar22}. Local regular solutions to problem \eqref{NS-IBVP}, as long as global solutions for small data, were provided in \cite{Mar-Pan} in the case of the Cauchy problem, while the linear problem in the half-space is discussed in \cite{AV}.}
\end{rem}
\begin{rem}
    {\rm It seems worth remarking that for initial data in the Fabes-Jones-Rivi\`ere class, that is $u_0\in L^q(\OO)$, $q>3$, a global weak solution can be obtained without many efforts. In fact, following the approach given in \cite{Kato}, we are able to obtain for this kind of initial datum a strong solution $u$, defined on some interval of time $(0,T_0)$. For the quoted solution, we find an instant of time $t_0\in (0,T_0)$ such that $u(t_0)\in W^{1,q}(\OO)\cap L^\infty(\OO)$. Therefore, we may take $u(t_0)$ as the initial datum for a global weak solution in the sense of Maremonti-Shimizu \cite{Mar-Shi-3D}.}
\end{rem}
\section{Preliminary results}\label{sec: prel}
We reproduce a result due to Calder\'{o}n in \cite{Cald} in the case $\OO=\R^3$. Anyway, we realize a decomposition of functions in general Lebesgue spaces without requiring a divergence-free condition to be preserved.
\begin{lemma}\label{Dec-Lp}
    {\sl Let $p,q,r\in [1,\infty)$ be such that $r<p<q$. For all $\rho>0$ and for all $f\in L^p(\Omega)$ there exist $f_1\in L^q(\Omega)$ and $f_2\in L^r(\Omega)$ such that $f=f_1+f_2$ and 
    \be\label{dec-est}
    \ba{rl}\displ
\dm f_1\dm_q&\hskip-0.2cm\leq \rho^{1-\frac{p}{q}}\dm f\dm_p^{\frac{p}{q}}\,,\VS \dm f_2\dm_r&\hskip-0.2cm\leq \rho^{1-\frac{p}{r}}\dm f\dm_p^{\frac{p}{r}}\,.
\ea
    \ee}
\end{lemma}
\bp
Let $\rho>0$. We write $f=f_1+f_2$ with 
\[
f_1(x)=\begin{cases}
    f(x)\,,\quad \text{if }|f(x)|\leq \rho\,,\\ 0\,,\quad \quad\,\,\text{else}\,.
\end{cases}
\]
Since $q-p>0$, it is clear that
\[
\dm f_1\dm_q^q=\int_\Omega|f_1(x)|^{q-p}|f_1(x)|^p\,dx\leq \rho^{q-p}\dm f\dm_p^p\,.
\]
Therefore, \eqref{dec-est}$_1$ holds. Analogously, since 
$r-p<0$ and $|f_2(x)|\geq \rho>0$ for all $x\in\OO$, we have
\[
\dm f_2\dm_r^r=\int_\Omega|f_2(x)|^{r-p}|f_2(x)|^p\,dx\leq \rho^{r-p}\dm f\dm_p^p\,,
\]
and \eqref{dec-est}$_2$ is also proved.
\ep
\begin{lemma}\label{le: Mar-AB}
    {\sl Let $y(t)\in C^{1}(t_0,\infty)$ such that $y(t)\geq 0$ and assume that
    \[
    \begin{cases}
        \dot y(t)\leq a(t)y(t)+K_1y^2(t)+K_2y^3(t)\,,\\ y(t_0)=y_0\,,\\  E=\int_{t_0}^\infty y(t)\,dt<\infty\,,
    \end{cases}
    \]
    for some $K_1,K_2\ge 0$. Assume further that there exists $K_3\ge 0$ such that $a(t)\leq K_3$ for all $t\ge t_0$, and that for some $\delta>0$
    \[
    E<\frac{\delta}{2K(1+\delta)^2}\,,\,\,\text{with }K:=\max\{K_i,\, i=1,\dots,3\}\,.
    \]
    Then $$y(t_0)\leq \frac\delta2\implies y(t)\leq \delta\,,\,\,\text{for all }t\ge t_0\,.$$ Moreover, if
    \[
    \int_{t_0}^\infty a(t)\,dt<\infty\,,
    \]
    there exists $c>0$ independent of $y$ such that
    \[
    y(t)\leq c(t-t_0)^{-1}\,,\,\,\text{for all }t>t_0\,.
    \]}
\end{lemma}
\bp
The lemma is proved in \cite{Mar-AB}
\ep
\begin{lemma}\label{le: CM}
    {\sl Assume that $D\subseteq \R^n$ is an exterior domain with compact boundary, having the cone property. Let $m\in \N_0$, $p\in [1,\infty]$, and $q\in[1,\infty)$. Let $g\in \widehat W^{m,p}(D)\cap L^q(D)$. Then, for $k\in\{0,1,\dots,m-1\}$, there exists $c>0$ independent of $g$ such that the following inequality holds:
    \be
\dm D^k g\|_r\leq c\dm D^m g\dm_p^a\dm g\dm_q^{1-a}\,,
    \ee
    where
    \[
    \frac1r=\frac{k}{n}+ a\bigg(\frac1p-\frac{m}{n}\bigg)+(1-a)\frac1q\,,
    \]
    with $a\in \big[\frac{k}{m},1\big]$ either if $p=1$ or if $p>1$ and $m-k-\frac{n}{p}\notin\N\cup\{0\}$, while $a\in \big[\frac{k}{m},1\big)$ if $p>1$ and $m-k-\frac{n}{p}\in\N\cup\{0\}$.}
\end{lemma}
\bp
See \cite{CM}.
\ep
\begin{lemma}\label{le: Sto-interp}
    {\sl Let $g\in J^{1,2}(\OO)\cap W^{2,2}(\OO)$. Then there exists a constant $c$, independent of $g$, such that
    \[
    \dm g\dm_\infty\leq c\dm\mathbb P_2\Delta g\dm_2^\frac12\dm \n g\dm_2^\frac12\,.
    \]}
\end{lemma}
\bp
See \cite{Mar-int}.
\ep
Let us consider the Stokes initial boundary value problem:
 \be\label{Sto-pb}
\begin{cases}
V_t-\Delta V+\n\pi_V=f\,,\,\,\quad\quad\text{in }(0,T)\times\OO,\\
\hskip1.83cm\n\cdot V=0\,,\,\,\quad\quad\text{in }(0,T)\times\OO\,,\\
\hskip1.6cmV(t,x)=0\,,\,\,\,\,\,\,\quad\,\,\text{in }(0,T)\times\partial\OO\,,\\
\hskip1.56cmV(0,x)=V_0(x)\,,\,\,\,\text{in }\Omega\,.
\end{cases}
\ee
\begin{lemma}\label{le: Stokes-Jp}
    {\sl Let  $\OO$ be an exterior domain with a $C^2$  boundary $\partial\OO$, $f=0$ and $V_0\in J^q(\OO)$, for some $q\in (1,\infty)$. Then, there exists a unique solution $(V,\pi_V)$ to problem \eqref{Sto-pb} such that, for all $T>0$, 
    \[
    \ba{c}
    \eta>0\,,\,\, V\in C([0,T);J^q(\OO))\cap L^q(\eta,T;J^{1,q}(\OO)\cap W^{2,q}(\OO))\,,\VS \n\pi_V,V_t\in L^q(\eta,T;L^q(\OO))\,.
    \ea
    \]
    Moreover, 
    \be \label{Sem-prop}
    \ba{rl}\displ
    \dm V(t)\dm_r&\hskip-0.2cm\leq c\dm V_0\dm_qt^{-\mu}\,,\quad \mu=\frac32\big(\frac1q-\frac1r\big)\,,\,\,t>0\,,\,r\in [q,\infty];\VS \dm \n V(t)\dm_r&\hskip-0.2cm\leq c\dm V_0\dm_qt^{-\mu_1}\,,\,\,\,\, \mu_1=\begin{cases}
        \frac12+\mu \,,\,\,\text{if }t\in(0,1)\,,\,r\in [q,\infty)\,,\\ \frac12+\mu\,,\,\,\text{if }t>0\,\,\text{and }r\in [q,3]\,,\\ \frac{3}{2q}\,,\quad\,\,\,\,\,\text{if }t\ge1\,\,\text{and }r> 3\,;
    \end{cases}\VS
    \dm V_t(t)\dm_r&\hskip-0.2cm\leq c\dm V_0\dm_qt^{-\mu_2}\,,\,\,\,\,\mu_2=1+\mu\,,\,\,t>0\,,\,r\in [q,\infty]\,;
    \ea
    \ee
    where the constant $c$ is independent of $V_0$.}
\end{lemma}
\bp
The Lemma is proved in \cite{Mar-L1}  by following analogous arguments to those employed in \cite{Mar-Sol}.
\ep
\begin{rem}\label{rem: Cauchy+semisp-Semi}
    {\rm In the case $\OO=\R^3$ or $\OO=\R^3_+$ the statement is clearly true and estimate \eqref{Sem-prop}$_2$ holds with $\mu_1=\frac12+\mu$. In the case of the Cauchy problem, this is an immediate consequence of the commutation between the Helmholtz projection and the Laplace operator, in the half-space case the result can be found in \cite{Ukai}. However, in the case of $\OO$ exterior domain estimate \eqref{Sem-prop}$_2$ is sharp, as proved in \cite{Mar-Sol}.}
\end{rem}
\begin{lemma}\label{le: pi-2}
    {\sl Let $f=0$ and $V_0\in J^2(\OO)$. Let $(V,\pi_V)$ be the solution to problem \eqref{Sto-pb} ensured by Lemma\,\ref{le: Stokes-Jp}. For all $\lambda\in (0,\frac12]$ and all $R_0>0$, setting  $\beta=\frac{1-2\lambda}{3}$, the pressure field $\pi_V$ satisfies 
    \[
    \dm \pi_V(t)\dm_{L^2(\Omega_{R_0})}\leq c(T)\dm V_0\dm_2[t^{-1+\frac\beta2}+t^{-1+\frac\beta2}]\,,\,\,\text{for all }t\in (0,T)\,.
    \]}
\end{lemma}
\bp
The proof is provided in \cite{CGM}.
\ep
\begin{coro}\label{cor: pi-2-extball}
    {\sl Under the same assumptions of Lemma\,\ref{le: pi-2}, for all $r>3$ the pressure field $\pi_V$ satisfies
    \[
    \dm \pi_V(t)\dm_{L^r(\Omega^{R_0})}\leq c(T)\dm V_0\dm_2[t^{-1+\frac\beta2}+t^{-1+\frac\beta2}]\,,\,\,\text{for all }t\in (0,T)\,.
    \]}
\end{coro}
\bp
The result is a consequence of Lemma\,\ref{le: pi-2}, as shown in \cite{CGM}. We briefly recall the key points here. We consider a smooth cut-off function $\varphi_{R_0}$, hence $\varphi_{R_0}(x)=0$ if $|x|\leq \frac{R_0}{2}$ and $\varphi_{R_0}(x)=1$ if $|x|\geq R_0$. We set $\zeta_{R_0}(x):=1-\varphi_{R_0}(x)$. Extending by zero the pressure field $\pi_V$ to the whole space, we consider, for all $t>0$, the Poisson equation
\[
\Delta(\zeta_{R_0}\pi_V)=2\n\zeta_{R_0}\cdot\n \pi_V+\Delta\zeta_{R_0}\pi_V\,,\,\,\text{for all }x\in\R^3\,.
\]
Denoting by $\mathcal E $ the fundamental solution of the Laplace equation, we have
\[
\zeta_{R_0}(x)\pi_V(x)=-\int_{\R^3}\mathcal E(x-y)\Delta \zeta_{R_0}(y)\pi_V(y)\,dy-2\int_{\R^3}\n\mathcal E(x-y)\n\zeta_{R_0}(y)\pi_V(y)\,dy\,.
\]
Since, for all $i\in\N$, $\text{supp}(\n^i\zeta_{R_0})\subset \{x\,:\, \frac{R_0}{2}<|x|<R_0\}$, the result follows from Lemma\,\ref{le: pi-2}.
\ep
\begin{lemma}\label{le: Sto-f}
    {\sl Let $V_0=0$ and $f\in L^q((0,T)\times\Omega)$, for some $q\in [2,\infty)$. Then there exists a unique solution $(V,\pi_V)$ to problem \eqref{Sto-pb} such that
    \[
    \ba{c}\displ
    V\in C([0,T);J^q(\OO))\cap L^q(0,T;J^{1,q}(\OO))\,,\VS \int_0^T\big[\dm V_t(\tau)\dm_q^q+\dm  D^2V(\tau)\dm_q^q+\dm \n\pi_V(\tau)\dm_q^q\big]\,d\tau\leq c(T)\int_0^T\dm f(\tau)\dm_q^q\,d\tau\,,
    \ea
    \]
    where, for all $\vep>0$, $c(T):=C(1+T^{1+\vep q-\frac{3}{2q}})$.}
\end{lemma}
\bp
See \cite{Mar-Sol}.
\ep
\begin{lemma}\label{le: Hardy-Ineq}
    {\sl Assume that for some $q\in (1,3)$ $\n f\in L^q(\OO)$. Then there exists a constant $f_0$ and  there exists $c>0$ independent of $f$ such that
    \[
    \dm f-f_0\dm_{\frac{3q}{3-q}}\leq c\dm\n f\dm_q\,.
    \]}
\end{lemma}
\bp
The lemma is proved in \cite{Ga:Book}.
\ep
\begin{lemma}\label{le: GWP-L3}
    {\sl There exists $\overline\xi_0>0$ such that each $v_0\in J^3(\Omega)$ with $\dm v_0\dm_3<\overline\xi_0$ corresponds to a unique solution $(v,\pi_v)$ to problem \eqref{NS-aux-intro} on $(0,\infty)$ such that, for all $T>\eta>0$,
    \begin{equation} \label{reg-sol-prel}
        \ba{c}
v\in C([0,T);J^3(\Omega))\cap L^\infty(\eta,T;J^{1,3}(\Omega)\cap W^{2,3}(\Omega))\,,\VS v_t,\n\pi_v\in L^\infty(\eta,T;L^3(\Omega))\,,\VS
t^{\frac32(\frac13-\frac1q)}v\in C([0,T);L^q(\Omega))\,,\,\,\text{for all }3\leq q\leq \infty\,,\VS
t^\frac12 \n v\in C([0,T);L^3(\Omega))\,, \VS
t^{\frac32(\frac13-\frac1q)+\frac12}\n v\in C([0,1);L^q(\OO))\,,\,\,\text{for all }3< q<\infty\,,\,\text{if }\OO\,\text{is an exterior domain}\,,\VS
 t^{\frac32(\frac13-\frac1q)+\frac12}\n v\in C([0,T);L^q(\Omega))\,,\,\,\text{for all }3<q<\infty\,, \,\,\text{if }\Omega=\R^3\,\text{or }\R^3_+\,.
\ea
    \end{equation}
    In particular, for suitable constants $h_1$, $h_2$ independent of $v_0$, and a constant $C=C(v_0)$ we have
    \begin{equation}\label{ASYMP-SOLGLO}
        \ba{rl}
        \dm v(t)\dm_q&\hskip-0.2cm\leq \frac{h_1\dm v_0\dm_3}{1+(1-h_2\dm v_0\dm_3)^\frac12}t^{-\frac32(\frac{1}{3}-\frac1q)}\,,\,\,\text{for all }q\in [3,\infty],\,\text{and }t>0\,,\VS
        \lim_{t\to\infty}\dm v(t)\dm_3&\hskip-0.2cm=0\,,\VS
        \dm v_t(t)\dm_3&\hskip-0.2cm\leq C\frac{h_1\dm v_0\dm_3}{1+(1-h_2\dm v_0\dm_3)^\frac12}t^{-1}\,,\,\,\text{for all }t>0\,.
        \ea
    \end{equation}
    Moreover, if $v_0\in J^3(\OO)\cap L^r(\OO)$, with $r\in (\frac32,3)$, then
    \begin{equation}\label{IDLrL3}
    \ba{rl}
        \dm v(t)\dm_q&\hskip-0.2cm\leq c\dm v_0\dm_rt^{-\frac32(\frac1r-\frac1q)}\,,\,\,\text{for all }t>0\,\,\text{and } q\in [r,\infty] \,,\VS \dm\n v(t)\dm_q&\hskip-0.2cm\leq \begin{cases} c\dm v_0\dm_rt^{-\frac32(\frac1r-\frac1q)-\frac12}\,,\,\,\text{for all }t>0\,\text{and }q\in [r,\infty)\,,\,\,\text{if }\OO=\R^3\,\text{or }\OO=\R^3_+,\\ c\dm v_0\dm_rt^{-\frac32(\frac1r-\frac1q)-\frac12}\,,\,\,\text{for all }t>0\,\text{and }q\in [r,3]\,,\,\,\text{if }\OO\,\text{is an exterior domain}.\end{cases}
        \ea
    \end{equation}}
\end{lemma}
\bp
The existence and uniqueness of a global regular solution, as long as properties \eqref{reg-sol-prel}$_{3,4,6}$ and \eqref{IDLrL3}, are furnished by Kato in \cite{Kato} in the context of the Cauchy problem. For the existence and regularity of solutions in the case of a bounded domain or $\OO=\R^3$ see also \cite{Giga86}, for the exterior problem see \cite{Iwa,Mar-DCDS}. The asymptotic properties \eqref{ASYMP-SOLGLO} are proved by Maremonti in \cite{Mar-DCDS}. Then we only need to prove the validity of properties \eqref{reg-sol-prel}$_{4,5}$ and \eqref{IDLrL3}$_{1,2}$ in the case of the initial boundary value problem. However, the arguments employed in \cite[Section 3]{Kato} still work in the case of a domain $\OO$ with a sufficiently smooth boundary. In fact, the proof furnished in the quoted reference requires both a smallness condition on the initial datum and the validity of the semigroup properties \eqref{Sem-prop}$_{1,2}$ for the Stokes operator.
As in our assumptions on the domain $\OO$ estimates \eqref{Sem-prop}$_{1,2}$ are satisfied, we get \eqref{reg-sol-prel}$_{4,5}$ and \eqref{IDLrL3}$_{2}$ as done in \cite{Kato}. In order to get \eqref{IDLrL3}$_1$, it is sufficient to employ a duality argument on the approximating sequence of the solution, as done in \cite{Mar-DCDS}. For further details we refer the reader to the appendix.
Concerning the half-space case, the result is proved by Kozono in \cite{Koz}.
\ep
\begin{rem} \label{rem: Cauchy VS IBVP}
    {\rm Actually, in the case of the Cauchy problem and the half-space problem, properties \eqref{IDLrL3} hold for all $r\in (1,3)$. This is due to the fact that $\mu_1=\frac12+\mu$, for all $t>0$ and $r\in [q,\infty)$, in \eqref{Sem-prop}$_2$. However, in the case of the exterior problem, properties \eqref{Sem-prop}$_{1,2}$ are sharp, as proved in \cite{Mar-Sol}. The different semigroup properties in the case of the exterior problem do not allow one to get \eqref{IDLrL3} for $r\in (1,\frac32)$. However, in order to reach our purposes, the requirement $r\in (\frac32,3)$ is  not restrictive, as the initial datum we are considering for the Navier-Stokes problem belongs to $L^p(\OO)$, $p\in (2,3)$. }
\end{rem}
As an immediate consequence of Lemma\,\ref{le: GWP-L3}, we get the following
\begin{coro}\label{cor: LTB-L3}
    {\sl Let $(v,\pi_v)$ the solution to problem \eqref{NS-aux-intro} obtained in Lemma\,\ref{le: GWP-L3}. Then
    \[
    \dm \n v(t)\dm_3\leq ct^{-\frac12}\dm v_0\dm_3\,,\,\,\text{for all }t>0\,.
    \]
    Moreover, if $v_0\in J^3(\OO)\cap L^r(\OO)$, for some $r\in (\frac32,3)$, then
    \[
     \dm \n v(t)\dm_3\leq ct^{-\frac{3}{2r}}\dm v_0\dm_r\,,\,\,\text{for all }t>0\,.
    \]}
\end{coro}
\begin{lemma}\label{le: EST-NONLIN}
    {\sl Let $w\in L^\infty(0,T;J^2(\Omega))\cap L^2(0,T;J^{1,2}(\Omega))$ with a bound $C$ independent of $T$. Let $q\in [\frac65, \frac{10}{7})$ and $v$ be the solution to problem \eqref{NS-aux-intro} given in Lemma\,\ref{le: GWP-L3}.
    Then, the following estimates hold:
    \begin{equation}\label{EST-NONLIN}
        \ba{c}\displ
        \int_0^T\dm w\cdot\n w\|_\frac54^\frac54\,d\tau\leq cC^\frac52\,,\quad \int_0^T\dm v\cdot\n w\dm_q^q\,d\tau\leq cC^q\dm v_0\dm_3^qT^\frac{10-7q}{4}\,,\VS \int_0^T\dm w\cdot\n v\dm_q^q\,d\tau\leq cC^q\dm v_0\dm_3^q\max\{T^\frac{10-7q}{4},1\}\,.
        \ea
    \end{equation}}
\end{lemma}
\bp
Property \eqref{EST-NONLIN}$_1$ was already proved by Ladyzhenskaya in \cite{Lad}. However, we propose a short proof here.
Employing Lemma\,\ref{le: CM}, we deduce that $\dm w\|_\frac{10}{3}\leq c\dm \n w\dm_2^\frac35\dm w\dm_2^\frac25$. By virtue of H\"{o}lder's inequality, we find that
\[
\int_0^T \dm w\cdot\n w\dm_\frac54^\frac54\,d\tau\leq \bigg[\int_0^T \dm w\dm_\frac{10}{3}^\frac{10}{3}\,d\tau\bigg]^{\frac38}\bigg[\int_0^T \|\n w\dm_2^2\,d\tau\bigg]^\frac58\leq cC^\frac52\,.
\]
Analogously, by virtue of \eqref{reg-sol}$_3$, we get
\[
\int_0^T \dm v\cdot\n w\dm_q^q\,d\tau\leq \bigg[\int_0^T \dm v\dm_\frac{2q}{2-q}^\frac{2q}{2-q}\,d\tau\bigg]^{\frac{2-q}{2}}\bigg[\int_0^T \|\n w\dm_2^2\,d\tau\bigg]^\frac{q}{2}\leq cC^q\dm v_0\dm_3^qT^\frac{10-7q}{4}\,,
\]
and, employing also \eqref{reg-sol}$_{4,5,6}$, we obtain
\[
\int_0^T \dm w\cdot\n v\dm_q^q\,d\tau\leq \int_0^T \dm w\dm_2^q\dm\n v\dm_\frac{2q}{2-q}^q\,d\tau\leq cC^q\dm v_0\dm_3^q\max\{1, T^\frac{10-7q}{4}\}\,.
\]
The lemma is completely proved.
\ep
\section{An auxiliary Navier-Stokes problem}\label{sec: auxprob}
In this section, we develop an analytic study of problem \eqref{NS-aux-intro}.
\begin{tho}\label{thm: L3-weakdata}
    {\sl There exists $\overline\xi_0>0$ such that for all $v_0\in L^3(\Omega)$ with $\dm \mathbb P_3v_0\dm_3<\overline\xi_0$, there exists a unique solution $(v,\pi_v)$ to problem \eqref{NS-aux-intro} on $(0,\infty)$ such that, for all $T>\eta>0$,
    \begin{equation}\label{reg-sol}
        \ba{c}
v\in C([0,T);J^3(\Omega))\cap L^\infty(\eta,T;J^{1,3}(\Omega)\cap W^{2,3}(\Omega))\,,\VS v_t,\n\pi_v\in L^\infty(\eta,T;L^3(\Omega))\,,\VS t^{\frac32(\frac13-\frac1q)}v\in C([0,T);L^q(\Omega))\,,\,\,\text{for all }3\leq q\leq \infty\,,\VS  t^\frac12\n v \in C([0,T);L^3(\OO))\,,\VS
t^{\frac32(\frac13-\frac1q)+\frac12}\n v\in C([0,1);L^q(\OO))\,,\,\,\text{for all }3< q<\infty\,,\,\text{if }\OO\,\text{is an exterior domain}\,,\VS
t^{\frac32(\frac13-\frac1q)+\frac12}\n v\in C([0,T);L^q(\Omega))\,,\,\,\text{for all }3< q<\infty\,, \,\,\text{if }\Omega=\R^3\,\text{or }\R^3_+\,.
        \ea
    \end{equation}
    In particular, for suitable constants $h_1$, $h_2$ independent of $v_0$, and a constant $C=C(v_0)$ we have
    \begin{equation}\label{ASYMP-SOLGLO-AUX}
        \ba{rl}
        \dm v(t)\dm_q&\hskip-0.2cm\leq \frac{h_1\dm \mathbb P_3v_0\dm_3}{1+(1-h_2\dm \mathbb P_3v_0\dm_3)^\frac12}t^{-\frac32(\frac{1}{3}-\frac1q)}\,,\,\,\text{for all }q\in [3,\infty],\,\text{and }t>0\,,\VS
        \lim_{t\to\infty}\dm v(t)\dm_3&\hskip-0.2cm=0\,,\VS
        \dm v_t(t)\dm_3&\hskip-0.2cm\leq C\frac{h_1\dm \mathbb P_3v_0\dm_3}{1+(1-h_2\dm \mathbb P_3v_0\dm_3)^\frac12}t^{-1}\,,\,\,\text{for all }t>0\,.
        \ea
    \end{equation} 
    Moreover, if $v_0\in L^3(\OO)\cap L^r(\OO)$, with $r\in (\frac32,3)$, then
    \begin{equation}\label{IDLrL3-aux}
    \ba{rl}
        \dm v(t)\dm_q&\hskip-0.2cm\leq c\dm v_0\dm_rt^{-\frac32(\frac1r-\frac1q)}\,,\,\,\text{for all }t>0\,\text{and }q\in[r,\infty]\,,\VS \dm\n v(t)\dm_q&\hskip-0.2cm\leq \begin{cases} c\dm v_0\dm_rt^{-\frac32(\frac1r-\frac1q)-\frac12}\,,\,\,\text{for all }t>0\,\text{and }q\in [r,\infty)\,,\,\,\text{if }\OO=\R^3\,\text{or }\OO=\R^3_+,\\ c\dm v_0\dm_rt^{-\frac32(\frac1r-\frac1q)-\frac12}\,,\,\,\text{for all }t>0\,\text{and }q\in [r,3]\,,\,\,\text{if }\OO\,\text{is an exterior domain}.\end{cases}
        \ea
    \end{equation}}
\end{tho}
\bp
We consider the unique global solution $(\overline v,\pi_{\overline v})$ to the related problem
\be\label{PROJ-Syst}
\begin{cases}
    \overline v_t-\Delta\overline v+\overline v\cdot\n\overline v+\n\pi_{\overline v}=0\,,\quad\,\,\,\,\,\,\,\,\,\,\,\quad \text{in }(0,T)\times\Omega\,,\\
   \hskip3.16cm \n\cdot\overline v=0\,,\quad\,\,\,\,\,\,\quad\,\,\,\,\,\text{in }(0,T)\times\Omega\,,\\
   \hskip2.89cm \overline v(t,x)=0\,,\quad \,\,\,\,\,\,\quad\quad\text{in }(0,T)\times\partial\Omega\,,\\
   \hskip2.84cm \overline v(0,x)=\mathbb P_3v_0\,,\quad\quad\,\,\,\text{in }\Omega\,,
\end{cases}
\ee
ensured by Lemma\,\ref{le: GWP-L3} thanks to the assumption $\dm \mathbb P_3 v_0\dm_3<\overline\xi_0$. By virtue of the Helmholtz decomposition, we have
\[
v_0(x)=\mathbb P_3v_0(x)+\n h\,.
\]
We recall that, by virtue of Lemma\,\ref{le: GWP-L3}, we have
\[
\lim_{t\to 0}\dm \overline v(t)-\mathbb P_3v_0\dm_3=0\,.
\]
Then, it follows that
\[
\lim_{t\to 0}(\overline v(t)-\mathbb P_3 v_0,\varphi)=0\,.
\]
Since
\[
(\n h,\varphi)=0\,,\,\,\text{for all }\varphi\in\mathscr C_0(\Omega)\,,\]
we deduce
\[
\lim_{t\to 0}(\overline v(t)- v_0,\varphi)=0\,,\,\,\text{for all }\varphi\in\mathscr C_0(\OO)\,.
\]
We conclude that $(\overline v,\pi_{\overline v})$ is also a global solution to problem \eqref{NS-aux-intro} and we refer to it as $(v,\pi_v)$. By virtue of Lemma\,\ref{le: GWP-L3}, it is clear that $(v,\pi_v)$ satisfies properties \eqref{ASYMP-SOLGLO-AUX}. It is also clear that $(v,\pi_v)$ satisfies \eqref{reg-sol} and, since the Helmholtz decomposition is well-defined in $L^r(\OO)$, with $r\in (\frac32,3)$, then if $v_0\in L^3(\OO)\cap L^r(\OO)$ we also have $\mathbb P_3 v_0\in J^3(\OO)\cap L^r(\OO)$. Therefore, estimates \eqref{IDLrL3-aux} hold if $v_0\in L^3(\OO)\cap L^r(\OO)$, $r\in (\frac32,3)$. We verify the uniqueness of the solution $(v,\pi_v)$. Actually, the proof is analogous to that of the uniqueness for problem \eqref{PROJ-Syst}, that follows an idea of Foias \cite{Foi}. We sketch it for the sake of completeness. We consider two solutions $(v_1,\pi_{v_1})$ and $(v_2,\pi_{v_2})$ to problem \eqref{NS-aux-intro} belonging to the class \eqref{reg-sol} and we set $$\widehat v=v_1-v_2\,,\quad \pi_{\widehat v}=\pi_{v_1}-\pi_{v_2}\,.$$
The functions $(\widehat v,\pi_{\widehat v})$ are solutions to the following problem:
\begin{equation}\label{DIFF-SYST}
    \begin{cases}
       \widehat v_t-\Delta\widehat v+\widehat v\cdot\n v_1+v_2\cdot\n\widehat v+\n\pi_{\widehat v}=0\,,\,\, \text{in }(0,T)\times\Omega\,,\\
   \hskip3.16cm \n\cdot\widehat v=0\,,\quad\,\,\,\,\,\,\quad\,\,\,\,\,\quad\text{in }(0,T)\times\Omega\,,\\
   \hskip2.89cm\widehat v(t,x)=0\,,\quad\quad \,\,\,\,\,\,\quad\quad\text{in }(0,T)\times\partial\Omega\,,\\
   \hskip2.14cm (\widehat v(0,x),\varphi)=0\,,\,\,\text{for all }\varphi\in\mathscr C_0(\Omega)\,.
    \end{cases}
\end{equation}
We consider $(z,\pi_z)$ solution to the Stokes problem with initial datum $z_0\in \mathscr C_0(\Omega)$. Let $t>0$. We set $\widehat z(\tau ,x)= z(t-\tau,x)$, $\tau \in (0,t)$. We multiply equation \eqref{DIFF-SYST}$_1$ by $\widehat z$ and integrate over $(0,t)\times \Omega$. We get
\be \label{dual-uniq}
(\widehat v(t),z_0)=(\widehat v(0),z(t))+\int_0^t \big[(\widehat v\cdot \n \widehat z, v_1)+(v_2\cdot\n\widehat z,\widehat v)\big]\,d\tau \,.
\ee
By virtue of condition \eqref{DIFF-SYST}$_4$, we have $(\widehat v(0),z(t))=0$. Moreover, we have
\[
\ba{l}\displ
\Bigg|\int_0^t \big[(\widehat v\cdot \n \widehat z, v_1)+(v_2\cdot\n\widehat z,\widehat v)\big]\,d\tau\Bigg|\leq \int_0^t\big|(\widehat v\cdot \n \widehat z, v_1)+(v_2\cdot\n\widehat z,\widehat v)\big|\,d\tau\VS\hskip5cm\leq \int_0^t\dm\widehat v(\tau)\dm_3\dm \n\widehat z(\tau)\dm_\frac32(\dm v_1(\tau)\dm_\infty+\dm v_2(\tau)\dm_\infty)\,d\tau\,.
\ea
\]
Set $$ k:=\int_0^t\tau^{-\frac12}(t-\tau)^{-\frac12}\,d\tau\,.$$
Employing the semigroup properties for $z$ and estimate \eqref{ASYMP-SOLGLO-AUX}$_1$ for $v_1$ and $v_2$, we obtain
\[
\ba{l}\displ
\int_0^t\dm\widehat v(\tau)\dm_3\dm \n\widehat z(\tau)\dm_\frac32(\dm v_1(\tau)\dm_\infty+\dm v_2(\tau)\dm_\infty)\,d\tau\VS\hskip1cm\leq c\dm z_0\dm_\frac32 \frac{h_1\dm\mathbb P_3 v_0\dm_3}{1+(1-h_2\dm\mathbb P_3 v_0\dm_3)^\frac12}\int_0^t\tau^{-\frac12}(t-\tau)^{-\frac12}\dm \widehat v(\tau)\dm_3\,d\tau\VS\hskip4cm\leq ck\dm z_0\dm_\frac32\frac{h_1\dm\mathbb P_3 v_0\dm_3}{1+(1-h_2\dm\mathbb P_3 v_0\dm_3)^\frac12}\sup_{(0,t)}\dm \widehat v(\tau)\dm_3\,.
\ea
\]
Since there is no lower bound for $\overline \xi_0$, we may assume that
\[
\overline \xi_0<(ck h_1)^{-1}\,,
\]
and, therefore,
\[
\overline c:=ck\frac{h_1\dm\mathbb P_3 v_0\dm_3}{1+(1-h_2\dm\mathbb P_3 v_0\dm_3)^\frac12}<1\,,
\]
Going back to \eqref{dual-uniq}, in view of the arbitrariness of $z_0$, we obtain
\[
(1-\overline c)\sup_{(0,t)}\dm \widehat v(\tau)\dm_3\leq 0\,,
\]
which in turn implies the uniqueness of the solution. The theorem is completely proved.
\ep
As a consequence of the previous theorem and Corollary\,\ref{cor: LTB-L3}, we have
\begin{coro}\label{cor: LTB-aux}
   {\sl Let $(v,\pi_v)$ the solution to problem \eqref{NS-aux-intro} obtained in Theorem\,\ref{thm: L3-weakdata}. Then
    \[
    \dm \n v(t)\dm_3\leq ct^{-\frac12}\dm v_0\dm_3\,,\,\,\text{for all }t>0\,.
    \]
    Moreover, if $v_0\in L^3(\OO)\cap L^r(\OO)$, for some $r\in (\frac32,3)$, then
    \[
     \dm \n v(t)\dm_3\leq ct^{-\frac{3}{2r}}\dm v_0\dm_r\,,\,\,\text{for all }t>0\,.
    \]}  
\end{coro}
\section{A perturbed Navier-Stokes problem}\label{sec: pertprob}
In this section, we perform an analytic study of problem \eqref{NS-pert-intro}. We recall that in the aforementioned problem $v$ denotes the solution to problem \eqref{NS-aux-intro}, ensured by Theorem\,\ref{thm: L3-weakdata}, and $w_0\in L^2(\Omega)$. Our goal is to construct a global weak solution to problem \eqref{NS-pert-intro} in the sense of Leray-Hopf. Further, we also discuss the partial regularity of the quoted solution, proving a structure theorem.
\begin{rem}\label{rem: costanti}
    {\rm We denote by $c_1$ the constant related to the embedding $\widehat W^{1,2}(\OO)\hookrightarrow L^6(\OO)$, namely 
\[
\dm f\dm_6\leq c_1\dm\n f\dm_2\,,\,\,\text{for all }f\in \widehat W^{1,2}(\OO)\,.
\]
As in Theorem\,\ref{thm: L3-weakdata} there is now lower bound for $\overline\xi_0$, we can take it so that
\[
\overline\xi_0<\min\{(c_1h_1)^{-1},h_2^{-1}\}\,,
\]
where $h_1$ and $h_2$ are given in \eqref{ASYMP-SOLGLO-AUX}$_1$.
Setting $\widetilde c:=c_1\frac{h_1}{1+(1-h_2\dm\mathbb P_3v_0\dm_3)^{\frac12}}$, we get 
\be \label{posdiss-energy}
1-\widetilde{c}\dm\mathbb P_3 v_0\dm_3>0\,.
\ee}
\end{rem}
We give the following definition for a weak solution.
\begin{defi}\label{def: WS-pert}
    {\sl A field $w:(0,\infty)\times\Omega\to\R^3$ is said a weak solution to problem \eqref{NS-pert-intro} if, for all $T>0$, the following properties hold:
    \begin{description}
        \item[(a)] $w\in L^\infty(0,T;L^2(\Omega))\cap L^2(0,T;W^{1,2}(\Omega))$;
        \item[(b)] for all $t,s\in (0,T)$, $w$ satisfies the following integral equation:
        \[
        \ba{l}\displ
        \int_s^t \big[(w,\phi_\tau)-(\n w,\n\phi)+((w+v)\cdot\n\phi, w)+(w\cdot\n\phi,v)\big]\,d\tau\VS\hskip0.5cm+ (w(s),\phi(s))=(w(t),\phi(t))\,,\,\,\text{for all }\phi\in C^1_0([0,T)\times \Omega),\,\text{with }\n\cdot\phi=0,\,\text{for all }t\in [0,T);
        \ea
        \]
        \item[(c)] $w$ satisfies the following inequality:
        \[
        \dm w(t)\dm_2^2+2(1-\widetilde c\dm\mathbb P_3v_0\dm_3)\int_s^t \dm \n w(\tau)\dm_2^2\,d\tau\leq \dm w(s)\dm_2^2\,,\,\,\text{for all }t\ge s\,\,\text{for }s=0\,\,\text{and for a. a. }s\ge 0\,;
        \]
        \item[(d)] $\displ\lim_{t\to 0^+}(w(t),\varphi)=(w_0,\varphi)$, for all $\varphi\in\mathscr C_0(\Omega)$.
    \end{description}}
\end{defi}
We also give the following notion of regular solution to problem \eqref{NS-pert-intro}
\begin{defi}\label{def: RSOL-pert}
    {\sl A weak solution $w$ to problem \eqref{NS-pert-intro} is said regular on $(0,T)\times\OO$ if 
    \[
    w\in L^2(0,T;W^{2,2}(\Omega))\,,\,\,w_t\in L^2(0,T;L^2(\Omega))
    \]
    and there exists a scalar field $\pi_w:(0,T)\times\OO\to\R$ such that
    \[
    \n\pi_w\in L^2(0,T;L^2(\OO))\,.
    \]}
\end{defi}
The main result of this section reads as follows.
\begin{tho}\label{thm: WSOL-pert}
    {\sl For all $w_0\in L^2(\OO)$ there exists a weak solution $w$ to problem \eqref{NS-pert-intro} in the sense of Definition\,\ref{def: WS-pert}. Moreover, for all $\eta>0$ there exist $\theta:= \eta^{-2}(1-\widetilde c\dm\mathbb P_3 v_0\dm_3)^{-1}\dm\mathbb P_2w_0\dm_2^4+1$, and a sequence of open time intervals $\{(\theta_l,T_l)\}_{l\in\N_0}$ such that 
    \[
    \ba{c}
    (\theta_0,T_0)=(\theta,\infty)\,, \,\,\text{and }(0,\infty)\setminus\underset{l\in\N_0}{\cup}(\theta_l,T_l)\,\,\text{has zero Lebesgue measure},\VS
    w\,\,\text{is regular on }(\theta_l,T_l)\times \OO\,,\,\,\text{for all }l\in\N_0\,.
    \ea
    \]}
\end{tho}
 The existence of a weak solution is proved following the idea of Leray \cite{Ler}, hence we study the following approximating problem. We set $\mathbb J_n[\cdot]=J_\frac1n[\cdot]$, where $J_\frac1n$ is a usual Friedrichs mollifier.
\begin{equation}\label{Moll-pb}
    \begin{cases}
        w^n_t-\Delta w^n+\mathbb J_n[w^n]\cdot\n w^n+w^n\cdot\n v+v\cdot\n w^n+\n\pi_{w^n}=0\,,\,\, \text{in }(0,T)\times\Omega\,,\\
   \hskip3.16cm \n\cdot w^n=0\,,\quad\,\,\,\,\,\,\quad\,\,\,\,\,\,\quad\qquad\qquad\,\,\,\,\qquad\quad\,\,\,\,\,\text{in }(0,T)\times\Omega\,,\\
   \hskip2.89cm w^n(t,x)=0\,,\quad\quad \,\,\,\,\,\,\quad\quad\qquad\qquad\,\,\,\,\,\,\qquad\quad\,\,\,\,\text{in }(0,T)\times\partial\Omega\,,\\
   \hskip2.84cm w^n(0,x)=\mathbb P_2w_0\,,\,\,\,\qquad\quad\,\,\,\,\,\qquad\quad\,\,\,\,\qquad\quad\,\,\,\,\,\,\text{in }\Omega\,.
    \end{cases}
\end{equation}
We recall again that $v$ is the regular solution to problem \eqref{NS-aux-intro}, ensured by Theorem\,\ref{thm: L3-weakdata}. \par
Problem \eqref{Moll-pb} can be studied employing the consolidated Faedo-Galerkin method in the way suggested by Prodi \cite{Pro} and Heywood \cite{Hey}. Since the procedure is well-understood, we omit it for the sake of brevity and we only perform the essential estimates on the Galerkin approximating sequences that allow to obtain a smooth solution to \eqref{Moll-pb}.
\subsection{Existence and uniqueness of a strong solution to \eqref{Moll-pb}}
In this subsection we establish the following result.
\begin{tho}\label{thm: WP-moll}
    {\sl For all $n\in\N$ there exists a unique solution $(w^n,\pi_{w^n})$ to problem \eqref{Moll-pb} such that 
    \[
    \ba{c}
    w^n\in C([0,T);J^2(\Omega))\cap L^2(0,T;J^{1,2}(\Omega))\,,\VS
    \text{for all }\delta>0,\, w^n_t,D^2w^n,\n\pi_{w^n}\in L^2(\delta,T;L^2(\Omega))\,.
    \ea
    \]}
\end{tho}
The proof of the theorem is based on \textit{a priori} estimates on the Galerkin approximating sequence $\{w^n_j\}_j$ related to problem \eqref{Moll-pb}. We recall that, given $R>\text{diam}(\R^3\setminus\OO)$, the $j$-th element of the Galerkin sequence $\{w^n_{j,R}\}_j$, on the invading domain $\OO_R$, is given by $w^n_{j,R}(t,x)=\displ \sum_{k=1}^j c^k_j(t)a_k(x)$, where $\{a_k\}_k$ is a special basis of $J^{1,2}(\OO_R)$, orthonormal in $J^2(\OO_R)$, with $\{a_k\}_k$ eigenfunctions of the Stokes operator $-\mathbb P_2\Delta$ and $\{\lambda_k\}_k$ the corresponding eigenvalues. The functions $c^k_j(t)$ satisfy the following system of ODEs:
\be\label{Gal-eq}
\begin{cases}
   ( \partial_t w^n_{j,R}-\mathbb P_2\Delta w^n_{j,R}-\mathbb J_n[w^n_{j,R}]\cdot\n w^n_{j,R}-v\cdot\n w^n_{j,R}-w^n_{j,R}\cdot\n v,a_k)=0,\\
    c^k_j(0)=(\mathbb P_2 w_0,a_k)\,, \,\,\, k=1,\dots,j\,.
\end{cases}
\ee
The estimates on the Galerkin approximating sequence are performed in the following lemmas, where, to avoid overloading the notation, we will omit the indexes $j$ and $R$. As a matter of fact, the estimates that will  be obtained will not depend on $(j,R)$. \par 
We first provide energy estimates.
\begin{lemma}\label{le: ENES}
    {\sl For all $n\in \N$, uniformly in $j\in\N$ and $R>\text{diam}(\R^3\setminus\OO)$, we get
    \begin{equation}
        \dm w^n(t)\dm_2^2+2(1-\widetilde c\dm \mathbb P_3v_0\dm_3)\int_0^t\dm \n w^n(t)\dm_2^2\,d\tau\leq\dm\mathbb P_2 w_0\dm_2^2\,, 
    \end{equation}
    where $\widetilde c$ is given in Remark\,\ref{rem: costanti}.}
\end{lemma}
\bp
Multiply equation \eqref{Gal-eq}$_1$ by $c^k_j$ and sum over the index $k=1,\dots,j$. Integrating over $(0,t)$, we get
\[
\dm w^n(t)\dm_2^2+2\int_0^t\dm \n w^n(\tau)\dm_2^2\,d\tau= \dm \mathbb P_2 w_0\dm_2^2+2\int_0^t(w^n\cdot\n v,w^n)\,d\tau\,.
\]
Integrating by parts, we obtain
\[
(w^n\cdot\n v,w^n)=-(w^n\cdot\n w^n,v)\,.
\]
Moreover, we have
\[
|(w^n\cdot\n w^n,v)|\leq c_1\dm v\dm_3\dm \n w^n\dm_2^2\leq \widetilde c \dm\mathbb P_3v_0\dm_3\dm \n w^n\dm_2^2\,.
\]
Therefore, we conclude the proof of the lemma.
\ep
In the following lemma we furnish estimates for second-order spatial derivatives and first-order time derivatives for the Galerkin sequence $\{w^n_{j,R}\}_j$.
\begin{lemma}\label{le: REG-APPR}
    {\sl Under the same assumptions of Lemma\,\ref{le: ENES}, for all $n\in\N$ we get, uniformly in $j\in \N$ and $R>\text{diam}(\R^3\setminus\OO)$, the following estimate:
    \begin{equation}\label{GRAD-TD}
    \ba{l}\displ
        t\dm \n w^n(t)\dm_2^2 + \int_0^t \tau \big[\dm \mathbb P_2\Delta w^n(\tau)\dm_2^2+\dm w^n_\tau(\tau)\dm_2^2\big]\,d\tau \VS\hskip1cm\leq c\int_0^t \tau\big[c(n)\dm w^n(\tau)\dm_2^2\dm\n w^n(\tau)\dm_2^2+(\dm \n v(\tau)\dm_3^2+\dm v(\tau)\dm_\infty^2)\dm\n w^n(\tau)\dm_2^2 \big]\,d\tau< \infty\,.
        \ea
    \end{equation}}
\end{lemma}
\bp
Firstly, we notice that the right-hand side of \eqref{GRAD-TD} is finite thanks to properties \eqref{reg-sol}$_{1,3,4,5}$ and Lemma\,\ref{le: ENES}. We now prove the validity of the estimate. We multiply equation \eqref{Gal-eq}$_1$ first by $-\lambda_k c^k_j$, then by $\tau$, sum over $k=1,\dots,j$ and integrate over $(0,t)$. We get
\be \label{iden-grad}
\ba{l}\displ
t\dm \n w^m(t)\dm_2^2+2\int_0^t\tau\dm \mathbb P_2\Delta w^m(\tau)\dm_2^2\,d\tau\VS\hskip1cm=2\int_0^t\tau\big[(\mathbb J_n[w^n]\cdot\n w^n,\mathbb P_2\Delta w^n)+(w^n\cdot\n v,\mathbb P_2\Delta w^n)+(v\cdot\n w^n,\mathbb P_2\Delta w^n)\big]\,d\tau\,.
\ea
\ee
We estimate the right-hand side of the previous identity employing elementary properties of mollifiers, H\"{o}lder and Young inequalities. We have
\[
\ba{rl}
|(\mathbb J_n[w^n]\cdot\n w^n,\mathbb P_2\Delta w^n)|&\hskip-0.2cm\leq c(n)\dm w^n\dm_2\dm \n w^n\dm_2\dm \mathbb P_2\Delta w^n\dm_2\leq \frac16\dm \mathbb P_2\Delta w^n\dm_2^2+c(n)\dm w^n\dm_2^2\dm \n w^n\dm_2^2\,,\VS
|(w^n\cdot\n v,\mathbb P_2\Delta w^n)|&\hskip-0.2cm\leq \dm \mathbb P_2\Delta w^n\dm_2\dm\n v\dm_3\dm w^n\dm_6\leq \frac16\dm \mathbb P_2\Delta w^n\dm_2^2+ c\dm\n v\dm_3^2\dm \n w^n\dm_2^2\,,\VS
|(v\cdot\n w^n,\mathbb P_2\Delta w^n)&\hskip-0.2cm \leq\dm v\dm_\infty\dm \n w^n\dm_2\dm\mathbb P_2 \Delta w^n\dm_2\leq\frac16\dm \mathbb P_2\Delta w^n\dm_2^2+c\dm v\dm_\infty^2\dm \n w^n\dm_2^2\,.
\ea
\]
Hence, we obtain
\[
\ba{l}\displ
\Bigg|\int_0^t\tau\big[(\mathbb J_n[w^n]\cdot\n w^n,\mathbb P_2\Delta w^n)+(w^n\cdot\n v,\mathbb P_2\Delta w^n)+(v\cdot\n w^n,\mathbb P_2\Delta w^n)\big]\,d\tau\Bigg|\VS\hskip0.2cm \leq \frac12\int_0^t \tau \dm\mathbb P_2\Delta w^n(\tau)\dm_2^2\,d\tau \!+\!\int_0^t\tau \big[c(n)\dm w^n(\tau)\dm_2^2\dm \n w^n(\tau)\dm_2^2 +c(\dm v(\tau)\dm_\infty^2+\dm \n v(\tau)\dm_3^2)\dm \n w^n(\tau)\dm_2^2 \big]d\tau\,.
\ea
\]
Therefore, going back to \eqref{iden-grad}, we deduce
\be\label{dis-grad}
 \ba{l}\displ
        t\dm \n w^n(t)\dm_2^2 + \int_0^t \tau \dm \mathbb P_2\Delta w^n(\tau)\dm_2^2\,d\tau \VS\hskip1cm\leq c\int_0^t \tau\big[c(n)\dm w^n(\tau)\dm_2^2\dm\n w^n(\tau)\dm_2^2+(\dm \n v(\tau)\dm_3^2+\dm v(\tau)\dm_\infty^2)\dm\n w^n(\tau)\dm_2^2 \big]\,d\tau\,.
        \ea
\ee
Analogously, multiplying equation \eqref{Gal-eq}$_1$ first by $\partial_\tau c^k_j$ and then by $\tau$, summing over $k=1,\dots,j$, and integrating over $(0,t)$, we get
\be\label{iden-TD}
\ba{l}\displ
t\dm \n w^m(t)\dm_2^2+2\int_0^t\tau\dm w^m_t(\tau)\dm_2^2\,d\tau\VS\hskip1cm=2\int_0^t\tau\big[(\mathbb J_n[w^n]\cdot\n w^n,w^n_t)+(w^n\cdot\n v, w^n_t)+(v\cdot\n w^n,w^n_t)\big]\,d\tau\,.
\ea
\ee
We estimate the right-hand side of \eqref{iden-TD} by analogous arguments to those employed for \eqref{iden-grad}. We finally deduce that
\be\label{dis-TD}
 \ba{l}\displ
        t\dm \n w^n(t)\dm_2^2 + \int_0^t \tau \dm  w^n_t(\tau)\dm_2^2\,d\tau \VS\hskip1cm\leq c\int_0^t \tau\big[c(n)\dm w^n(\tau)\dm_2^2\dm\n w^n(\tau)\dm_2^2+(\dm \n v(\tau)\dm_3^2+\dm v(\tau)\dm_\infty^2)\dm\n w^n(\tau)\dm_2^2 \big]\,d\tau\,.
        \ea
\ee
Combining estimates \eqref{dis-grad} and \eqref{dis-TD}, we conclude the proof.
\ep
\bp[Proof of Theorem\,\ref{thm: WP-moll}]
The existence of a solution $(w^n,\pi_{w^n})$ such that
\[
\ba{c}
w^n\in C([0,T);J^2(\OO))\cap L^2(0,T;J^{1,2}(\OO))\,,\VS w^n_t,D^2 w^n,\n \pi_{w^n}\in L^2(\delta,T;L^2(\OO))\,,\,\,\text{for all }\delta >0\,,
\ea
\]
can be achieved by following the procedure suggested in \cite{Hey, Pro} in view of Lemmas\,\ref{le: ENES}-\ref{le: REG-APPR}. In order to prove the uniqueness of the quoted solution, we assume that $(w^n_1,\pi_{w^n_1})$ and $(w^n_2,\pi_{w^n_2})$ are two solutions to problem \eqref{Moll-pb}. We set $(\widehat w^n,\widehat \pi^n):=(w^n_1-w^n_2,\pi_{w^n_1}-\pi_{w^n_2})$. It is readily seen that $(\widehat w^n,\widehat \pi^n)$ satisfies
\be \label{Diff-moll}
\hskip-0.15cm\begin{cases}
       \widehat w^n_t\!-\!\Delta\widehat w^n\!+\!\mathbb J_n[\widehat w^n]\cdot\n w^n_1\!+\!\mathbb J_n[w^n_2]\cdot \n\widehat w^n\!+\!\widehat w^n\cdot\n v\!+\!v\cdot\n\widehat w^n\!+\!\n\widehat\pi^n=0\,,\,\, \text{in }(0,T)\times\Omega\,,\\
   \hskip3.16cm \n\cdot\widehat w^n=0\,,\quad\,\,\,\,\,\,\quad\,\,\,\,\,\,\quad\qquad\qquad\,\,\,\,\qquad\quad\,\,\,\,\,\text{in }(0,T)\times\Omega\,,\\
   \hskip2.91cm\widehat w^n(t,x)=0\,,\quad\quad \,\,\,\,\,\,\quad\quad\qquad\qquad\,\,\,\,\,\,\qquad\quad\,\,\,\,\text{in }(0,T)\times\partial\Omega\,,\\
   \hskip2.86cm\widehat w^n(0,x)= 0\,,\,\,\,\qquad\quad\,\,\,\,\,\qquad\quad\,\,\,\,\qquad\quad\,\,\,\,\,\,\quad\,\,\,\,\text{in }\Omega\,.
    \end{cases}
\ee
We multiply equation \eqref{Diff-moll}$_1$ by $\widehat w^n$ and integrate over $\OO$. By virtue of the divergence theorem, we get
\[
\frac12\frac{d}{dt}\dm \widehat w^n(t)\dm_2^2+\dm \n \widehat w^n(t)\dm_2^2=-(\mathbb J_n[\widehat w^n(t)]\cdot\n w^n_1(t),\widehat w^n(t))-(\widehat w^n(t)\cdot\n v(t),\widehat w^n(t))\,,\,\text{for all }t>0\,.
\]
We estimate the right-hand side of the above identity. We have
\[
\ba{rl}
|(\mathbb J_n[\widehat w^n]\cdot\n w^n_1,\widehat w^n)|&\hskip-0.2cm\leq c(n)\dm \n w^n_1\dm_2\dm \widehat w^n\dm_2^2\,,\VS
|(\widehat w^n\cdot\n v,\widehat w^n)|=
|(\widehat w^n\cdot\n\widehat w^n,v)|&\hskip-0.2cm\leq \dm v\dm_3\dm \n \widehat w^n\dm_2\dm \widehat w^n\dm_6\leq \widetilde c\dm\mathbb P_3 v_0\dm_3\dm\n \widehat w^n\dm_2^2\,,
\ea
\]
where $\widetilde c$ is given in Remark\,\ref{rem: costanti}. Therefore, recalling the aforementioned remark and that, by virtue of Lemma\,\ref{le: ENES}, $\dm \n w^n_1\dm_2\in L^2(0,T)$, for all $T>0$, we conclude the proof of the theorem by invoking Gronwall's Lemma.
\ep
\begin{rem}
    {\rm It seems worth to emphasize that, differently from what is known for the $L^3$-theory for global weak solutions (see \cite{Mar-L3Lor,Ser-Sve}), with our technique the approximating sequence for the weak solution $w$ is uniquely determined.}
\end{rem}
\subsection{Existence of a weak solution to problem \eqref{NS-pert-intro}}
We want to prove the existence of a weak solution to problem \eqref{NS-pert-intro}. The weak solution will be obtained  through a limit process on the sequence $\{w^n\}_n$, constructed in the previous subsection. In order to reach our goal we need two preparatory lemmas. The former allows us to deduce some estimates concerning the pressure fields $\{\pi_{w^n}\}_n$, uniformly with respect to $n\in\N$. The latter shows that, for all $t>0$, the $L^2$ norm of $w^n(t)$ in the exterior domain $\Omega^R$ vanishes letting $R\to\infty$, uniformly in $n\in\N$. 
\begin{lemma}\label{le: Press-moll}
    {\sl Let $(w^n,\pi_{w^n})$ be the solution to problem \eqref{Moll-pb} established in Theorem\,\ref{thm: WP-moll} and let $\lambda\in (0,\frac12)$. Set $\beta=\frac{1-2\lambda}{3}$. Then, $\pi_{w^n}=\pi^n_1+\pi_2^n$, with
    \be \label{press-est-moll}
    \ba{c}\displ
    \int_0^T\dm \n \pi_{1}^n(\tau)\dm_\frac{5}{4}^\frac54\,d\tau\leq C_1(\dm v_0\dm_3,T)\,,\VS\dm\pi_2^n(t)\dm_{L^r(\OO^{R})}\leq c(T)\dm\mathbb P_2 w_0\dm_2[t^{-1+\frac\beta2}+t^{-1+\frac\beta4}]\,,\,\,\text{for all }R>0\,, \,\,r>3\,\, \,\,\text{and }t\in (0,T)\,.
    \ea
    \ee
    }
    \end{lemma}
    \bp
    We can write $(w^n,\pi_{w^n})=(w^n_1,\pi_1^n)+(w_2^n,\pi^n_2)$, where $(w^n_1,\pi^n_1)$ solves 
    \begin{equation}\label{Sto-w1}
        \begin{cases}
            \partial_tw^n_1-\Delta w^n_1=-\n\pi_1^n+f^n\,,\quad \text{in }(0,T)\times\OO\,,\\
            \hskip1cm\n\cdot w^n_1=0\,,\quad\qquad\quad\,\,\,\quad\text{in }(0,T)\times\OO\,,\\
            \hskip0.8cmw^n_1(t,x)=0\,,\quad\qquad\quad\,\,\,\quad\text{in }(0,T)\times\partial\OO\,,\\ \hskip0.75cmw^n_1(0,x)=0\,,\quad\qquad\quad\,\,\,\quad\text{in }\{0\}\times\OO\,\,,
        \end{cases}
    \end{equation}
    with
    \[
    f^n=-\mathbb J_n[w^n]\cdot\n w^n-v\cdot\n w^n-w^n\cdot\n v\,,
    \]
    while $(w^n_2,\pi^n_2)$ solves
    \begin{equation}\label{Sto-w2}
        \begin{cases}
            \partial_tw^n_2-\Delta w^n_2=-\n\pi_2^n\,,\,\,\quad\quad\quad \text{in }(0,T)\times\OO\,,\\
            \hskip1cm\n\cdot w^n_2=0\,,\quad\qquad\quad\,\,\,\quad\text{in }(0,T)\times\OO\,,\\
            \hskip0.8cmw^n_2(t,x)=0\,,\quad\qquad\quad\,\,\,\quad\text{in }(0,T)\times\partial\OO\,,\\ \hskip0.75cmw^n_2(0,x)=\mathbb P_2 w_0\,,\quad\,\,\quad\,\,\,\quad\text{in }\{0\}\times\OO\,\,.
        \end{cases}
    \end{equation}
    Therefore, estimate \eqref{press-est-moll}$_1$ is a consequence of Lemma\,\ref{le: EST-NONLIN} and Lemma\,\ref{le: Sto-f}, while estimate \eqref{press-est-moll}$_2$ follows from Lemma\,\ref{le: pi-2} and Corollary\,\ref{cor: pi-2-extball}. The lemma is proved.
    \ep
    \begin{rem}\label{rem: Hardy-press}
        {\rm We point out that, by virtue of Lemma\,\ref{le: Hardy-Ineq}, there exists a function $c_1(t)$ such that 
        \[
        \int_0^T \dm \pi^n_1(\tau)-c_1(\tau)\dm_\frac{15}{7}^\frac54\,d\tau\leq C_1(\dm v_0\dm_3,T)\,.
        \]
        From now on, with an abuse of notation, we will denote the field $\pi^n_1(t,x)-c_1(t)$ employing the concise notation $\pi^n_1(t,x)$. }
    \end{rem}
   
    \begin{lemma}\label{le: Conv-extball}
        {\sl Let $(w^n,\pi_{w^n})$ be the solution to problem \eqref{Moll-pb} established in Theorem\,\ref{thm: WP-moll}. Then, for all $t>0$, we have
        \be
\lim_{R\to\infty}\dm w^n(t)\dm_{L^2(\OO^R)}=0\,,\,\,\text{uniformly in }n\in\N\,.
        \ee}
    \end{lemma}
    \bp
Consider a smooth cut-off function $\varphi_R(x)$, hence $\varphi_R(x)=1$ if $|x|\leq R$ and $\varphi_R(x)=0$ if $|x|\geq 2R$. Moreover, $\text{supp}(\n^i\zeta_R)\subset\{x\in\R^3\,:\, R<|x|<2R\}$ and $|\n^i\varphi_R(x)|\leq cR^{-i}$, for all $x\in \R^3$, $i\in \N$. Set $\zeta_R(x):=1-\varphi_R(x)$. We multiply equation \eqref{Moll-pb}$_1$ by $\zeta_R^2w^n$. Integration by parts over $(0,t)\times\OO$ yields
\be\label{w-extsphere}
\ba{l}\displ
\frac12\dm \zeta_Rw^n(t)\dm_2^2\!+\!\int_0^t\dm \zeta_R^2\n w^n(\tau)\dm_2^2\,d\tau=\!\int_0^t (|w^n(\tau)|^2,\Delta\zeta_R^2)\,d\tau\VS\hskip1cm+\!\int_0^t ([\mathbb J_n[w^n]+v(\tau)]\cdot \n \zeta_R^2,|w^n(\tau)|^2)\,d\tau+\int_0^t(w^n(\tau)\cdot\n v(\tau),\zeta_R^2w^n(\tau))\,d\tau\VS\hskip5cm+\int_0^t(\pi_{w^n}(\tau)w^n(\tau),\n\zeta_R^2)\,d\tau=:\sum_{i=1}^4I_i(t,n,R)\,.
\ea
\ee
By virtue of Lemma\,\ref{le: ENES}, we get
\[
|I_1(t,n,R)|\leq c(t,\dm \mathbb P_2 w_0\dm_2^2)R^{-2}\,.
\]
Concerning $I_2$, employing H\"{o}lder's inequality we have
\[
|I_2(t,n,R)|\leq cR^{-1}\int_0^t\big[\dm w^n(\tau)\dm_3^3+\dm v\dm_\infty\dm w^n(\tau)\dm_2^2\big]\,d\tau\,.
\]
Since $\dm w^n\dm_3\leq \dm w^n\dm_2^\frac12\dm\n w^n\dm_2^\frac12$ and, by virtue of property \eqref{ASYMP-SOLGLO-AUX}$_1$, $v\in L^1(0,T;L^\infty(\OO))$, invoking Lemma\,\ref{le: ENES}, we deduce
\[
|I_2(t,n,R)|\leq c(t,\dm\mathbb P_2w_0\dm_2,\dm v_0\dm_3 )R^{-1}\,.
\]
Concerning $I_3$, integrating by parts we obtain
\[
I_3(t,n,R)=-\int_0^t\big[(w^n(\tau)\cdot\n w^n(\tau),\zeta_R^2v(\tau))+(w^n(\tau)\cdot (\n\zeta_R^2\otimes w^n(\tau)),v(\tau))\big]\,d\tau=:-(I_{3,1}+I_{3,2})\,.
\]
Employing H\"{o}lder's inequality and recalling Remark\,\ref{rem: costanti}, we get
\[
|I_{3,1}(t,n,R)|\leq \int_0^t \dm w^n(\tau)\dm_6\dm \zeta_R^2\n w^n(\tau)\dm_2\dm v(\tau)\dm_3\,d\tau\leq \widetilde c\int_0^t \dm\mathbb P_3 v_0\dm_3\dm \zeta_R^2\n w^n(\tau)\dm_2^2\,d\tau\,.
\]
Moreover, concerning $I_{3,2}$, employing H\"{o}lder's inequality, we find that
\[
|I_{3,2}(t,n,R)|\leq cR^{-1}\int_0^t\dm w^n(\tau)\dm_6\dm\n w^n(\tau)\dm_2\dm v(\tau)\dm_3\,d\tau.
\]
Invoking again Lemma\,\ref{le: ENES} and \eqref{ASYMP-SOLGLO-AUX}$_1$, we conclude that
\[
\lim_{R\to \infty}|I_{3,2}(t,n,R)|=0\,,\quad\text{uniformly in }n\in\N\,.
\]
Finally, we discuss the term $I_4$. We have
\[
I_4(t,n,R)=\int_0^t (\pi_{w^n}(\tau)\n\zeta_R^2,w^n(\tau))\,d\tau=\int_0^t ((\pi_{1}^n(\tau)+\pi^n_2(\tau))\n\zeta_R^2,w^n(\tau))\,d\tau=:I_{4,1}+I_{4,2}\,,
\]
where
\[
\ba{rl}
I_{4,1}(t,n,R)&\hskip-0.2cm:=\displ\int_0^t (\pi_{1}^n(\tau)\n\zeta_R^2,w^n(\tau))\,d\tau\,,\VS
I_{4,2}(t,n,R)&\hskip-0.2cm:=\displ\int_0^t (\pi_{2}^n(\tau)\n\zeta_R^2,w^n(\tau))\,d\tau\,.
\ea
\]
By virtue of H\"{o}lder's inequality, Lemma\,\ref{le: ENES} and the properties of cut-off functions, we get
\[
\ba{l}\displ
\Bigg|\int_0^t (\pi_{1}^n(\tau)\n\zeta_R^2,w^n(\tau))\,d\tau\Bigg|\displ\leq \int_0^t\dm \pi^n_1(\tau)\dm_\frac{15}{7}\dm \n \zeta_R\dm_{30}\dm w^n(\tau)\dm_2\,d\tau\VS\hskip1cm\leq cR^{-\frac{9}{10}}\int_0^t\dm \pi^n_1(\tau)\dm_\frac{15}{7}\dm w^n(\tau)\dm_2\,d\tau\leq c(\dm\mathbb P_2w_0\dm_2)R^{-\frac{9}{10}}\int_0^t\dm \pi^n_1(\tau)\dm_\frac{15}{7}\,d\tau\,.
\ea
\]
Recalling Remark\,\ref{rem: Hardy-press}, we conclude that
\[
\lim_{R\to \infty}|I_{4,1}(t,n,R)|=0\,,\quad\text{uniformly in }n\in\N\,.
\]
Finally, we consider $I_{4,2}$. Let $r\in (3,6)$. Since $\text{supp}(\n\zeta_R)\subset \{x\,:\,R<|x|<2R\}\subset \OO^R$, employing H\"{o}lder's inequality, Lemma\,\ref{le: Press-moll} and the properties of cut-off functions, we get
\[
|I_{4,2}(t,n,R)|\leq \int_0^t\dm w^n(\tau)\dm_2\dm \n \zeta_R\dm_{\frac{2r}{r-2}}\dm \pi^n_2(\tau)\dm_{L^r(\OO^R)}\,d\tau\leq c\dm\mathbb P_2w_0\dm_2^2R^{\frac{r-6}{2r}}\int_0^t[\tau^{-1+\frac\beta2}+\tau^{-1+\frac\beta4}]\,d\tau\,.
\]
We recall that from Lemma\,\ref{le: Press-moll} we have $\beta>0$. Since 
\[
\int_0^t[\tau^{-1+\frac\beta2}+\tau^{-1+\frac\beta4}]\,d\tau=\frac2\beta t^{\frac\beta2}+\frac4\beta t^{\frac\beta4}\,,
\]
and $r\in (3,6)$,
we conclude that
\[
\lim_{R\to \infty}|I_{4,2}(t,n,R)|=0\,,\quad\text{uniformly in }n\in\N\,.
\]
Going back to \eqref{w-extsphere}, we deduce
\[
\frac12\dm \zeta_Rw^n(t)\dm_2^2\!+(1-\widetilde c\dm\mathbb P_3 v_0\dm_3)\!\int_0^t\dm \zeta_R^2\n w^n(\tau)\dm_2^2\,d\tau\leq c \big(\sum_{i=1,2,4}|I_i(t,n,R)|+|I_{3,2}(t,n,R)|\big)\,.
\]
As we proved that the right-hand side of the above inequality tends to $0$, uniformly in $n\in\N$, as $R\to\infty$, we
conclude the proof of the Lemma.
    \ep
    We are now in a position to prove the main result of this subsection.
    \begin{tho}\label{thm: WSOL-per-exist}
        {\sl There exists a weak solution $w$ to problem \eqref{NS-pert-intro} in the sense of Definition\,\ref{def: WS-pert}.}
    \end{tho}
    \bp
We deduce our weak solution $w$ as the weak limit of the sequence $\{w^n\}_{n\in\N}$, established in Theorem\,\ref{thm: WP-moll}, in a suitable metric.
        By virtue of Lemma\,\ref{le: ENES}, we see that the sequence $\{w^n\}_n$ is bounded in $L^2(0,T;J^{1,2}(\OO))$. Hence, there exists a subsequence, still denoted by $\{w^n\}$, which admits a weak limit $w$ in this space.\par 
        Employing consolidated arguments, see \cite{Ga: LecNot, Lad, Pro}, one can prove that, for all $\psi\in J^2(\OO)$, the sequence
        \[
        \{(w^n(\tau),\psi)\}_{n\in\N}
        \]
        is equicontinuous and uniformly bounded. Hence, the Ascoli-Arzel\`a Theorem ensures that there exists a subsequence weakly converging to a function $h(t)\in C([0,T))$. Then, we deduce that 
        \[
        h(t)=(w(t),\psi)\,,\,\,\text{a. e. in }[0,T)\,.
        \]
         We now show that $w$ satisfies the strong energy inequality given in the item $\bf(c)$ of Definition\,\ref{def: WS-pert}. 
        It is readily seen that the sequence $\{w^n\}_n$ satisfies the energy relation
        \be \label{SEI-moll}
        \dm w^n(t)\dm_2^2+2(1-\widetilde c\dm\mathbb P_3v_0\dm_3)\int_s^t \dm\n w^n(\tau)\dm_2^2\,d\tau\leq \dm w^n(s)\dm_2^2\,, \,\,\text{for all }t>s\,\,\text{and a. a. }s\geq 0\,.
        \ee
        As a matter of fact, the above inequality is obtained multiplying equation \eqref{Moll-pb}$_1$ by $w^n$, integrating over $(s,t)\times\OO$, employing H\"{o}lder's inequality as long as estimate \eqref{ASYMP-SOLGLO-AUX}$_1$, and recalling Remark\,\ref{rem: costanti}.
        We show that for a. a. $s\in (0,t)$
        \[
        \lim_{n\to\infty} \dm w^n(s)\dm_2^2=\dm w(s)\dm_2^2\,.
        \]
       In fact, since $$\dm w^n(s)\dm_2^2= \dm w^n(s)\dm_{L^2(\OO_R)}^2+\dm w^n(s)\dm_{L^2(\OO^R)}^2\,, $$
        and, by virtue of Friedrich's Lemma and Lemma\,\ref{le: Conv-extball} we have $$\ba{c} \displ w^n\to w\,\,\text{strongly in }L^2(0,T;L^2(\OO_R))\,, \VS\lim_{R\to\infty}\dm w^n(s)\dm_{L^2(\OO^R)}=0\,,\,\,\text{uniformly in }n\in\N\,, 
        \ea$$  and our claim is achieved. 
        We take the minimum limit in \eqref{SEI-moll} and we get
        \[
        \dm w(t)\dm_2^2+2(1-\widetilde c\dm\mathbb P_3v_0\dm_3)\int_s^t\dm \n w(\tau)\dm_2\,d\tau\leq \dm w(s)\dm_2^2\,,\,\,\text{for all }t>0\,\text{and a. a. }s\in (0,t)\,.
       \]
        We conclude that $w$ satisfies item $\bf (c)$ of Definition\,\ref{def: WS-pert}. From the validity of item $\bf (c)$ we also get that item $\bf (a)$ is satisfied. The proof of items $\bf (d)$ and $\bf (b)$ can be achieved employing classical arguments, see, for instance, \cite{Lad}.
    \ep
    \subsection{Partial regularity of the weak solution}
    We now provide for the weak solution obtained in Theorem\,\ref{thm: WSOL-per-exist} a \textit{Th\'{e}or\`{e}me de Structure} in the spirit of the celebrated result of Leray for the Navier-Stokes problem. We follow the ideas introduced by P. Maremonti and the author of this manuscript in the recent paper \cite{Mar-Pal-NS}\footnote{The new proof of the \textit{Th\'{e}or\`{e}me de Structure} in the paper \cite{Mar-Pal-NS} is based on arguments that might be employed also in different context of parabolic problems, in this connection see \cite{Mar-Pal-FSI} for applications to the dynamic interaction between a rigid body and a viscous incompressible fluid. }. We first need to prove the following
    \begin{lemma}\label{le: SECDER-est}
    {\sl Let $\{w^n\}_{n\in\N}$ be the sequence established in Theorem\,\ref{thm: WP-moll}. Then, for all $\tau >0$ and $n\in\N$, the following inequality holds:
    \be \label{ineq-grad}
\frac{d}{dt}\dm \n w^n(\tau)\dm_2^2+\frac12\dm \mathbb P_2\Delta w^n(\tau)\dm_2^2+\dm w^n_\tau(\tau)\dm_2^2\leq c(\dm v(\tau)\dm_\infty^2+\dm \n v(\tau)\dm_3^2+\dm \n w^n(\tau)\dm_2^4)\dm \n w^n(\tau)\dm_2^2\,.
    \ee
    Moreover, for all $\vep\in (0,1)$ and $T>0$, 
    \be\label{est-DERSEC}
\int_0^T \dm \mathbb P_2\Delta w^n(t)\dm_2^{\frac{2}{3+\vep}}\,dt\leq C(T,\dm \mathbb P_2 w_0\dm_2)<\infty\,, \,\,\text{for all }n\in\N\,.
    \ee}
    \end{lemma}
    \bp
We apply the Helmholtz projection operator $\mathbb P_2$ to both sides of equations \eqref{Moll-pb}$_1$ and then evaluate the $L^2$ norm of both sides. For all $\tau>0$, we get
\be\label{DER-sec-step1}
\ba{l}\displ
\frac{d}{dt}\dm \n w^n(\tau)\dm_2^2+\dm \mathbb P_2\Delta w^n(\tau)\dm_2^2+\dm w^n_\tau(\tau)\dm_2^2=\dm \mathbb P_2[\mathbb J_n[w^n(\tau)]\cdot\n w^n(\tau)]\dm_2^2\VS\hskip4cm+\dm \mathbb P_2[v(\tau)\cdot\n w^n(\tau)]\dm_2^2 +\dm \mathbb P_2[w^n(\tau)\cdot\n v(\tau)]\dm_2^2\,.
\ea
\ee
We estimate the right-hand side of \eqref{DER-sec-step1}. Employing H\"{o}lder's inequality and recalling Lemma\,\ref{le: Sto-interp}, we have
\[
\ba{rl}
\dm \mathbb P_2[\mathbb J_n[w^n(\tau)]\cdot\n w^n(\tau)]\dm_2^2&\hskip-0.2cm\leq c\dm w^n(\tau)\dm_\infty^2\dm \n w^n(\tau)\dm_2^2\leq \frac12\dm\mathbb P_2\Delta w^n(\tau)\dm_2^2+c\dm \n w^n(\tau)\dm_2^6\,,\VS
\dm \mathbb P_2[v(\tau)\cdot\n w^n(\tau)]\dm_2^2&\hskip-0.2cm\leq c\dm v(\tau)\dm_\infty^2\dm\n w^n(\tau)\dm_2^2\,,\VS 
\dm \mathbb P_2[w^n(\tau)\cdot\n v(\tau)]\dm_2^2&\hskip-0.2cm\leq c\dm w^n(\tau)\dm_6^2\dm \n v(\tau)\dm_3^2\leq c\dm \n v(\tau)\dm_3^2\dm\n w^n(\tau)\dm_2^2\,.
\ea
\]
Therefore, we deduce the validity of \eqref{ineq-grad}.  Let $\alpha\in(0,\vep)$ and multiply both sides of \eqref{ineq-grad} by $\tau^\alpha (1+\dm \n w^n(\tau)\dm_2^2)^{-2}$. Employing \eqref{reg-sol}$_{1,3,4}$ for $v(\tau)$, $\tau>0$, we obtain
\[
\tau^\alpha\dm \mathbb P_2\Delta w^n(\tau)\dm_2^2(1+\dm \n w^n(\tau)\dm_2^2)^{-2}\leq a\tau^{\alpha-1}+c\tau^{\alpha}\dm \n w^n(\tau)\dm_2^2-\tau^\alpha \frac{d}{d\tau}(1+\dm\n w^n(\tau)\dm_2^2)^{-1}\,.
\]
We integrate both sides of the above inequality on $(0,T)$ and apply to the left-hand side the reverse H\"{o}lder's inequality with exponents $\frac{1}{3+\vep}$ and $-\frac{1}{2+\vep}$. We get
\be \label{Rev-Hold}
\ba{l}\displ
\bigg(\int_0^T\dm \mathbb P_2\Delta w^n(\tau)\dm_2^\frac{2}{3+\vep}\,d\tau\bigg)^{3+\vep}\bigg(\int_0^T\tau^{-\frac{\alpha}{2+\vep}}(1+\dm \n w^n(\tau)\dm_2^2)^{\frac{2}{2+\vep}}\,d\tau\bigg)^{-(2+\vep)}\VS\hskip1cm\leq \bigg[\frac{a}{\alpha}T^\alpha+ cT^\alpha\int_0^T\dm \n w^n(\tau)\dm_2^2\,d\tau-\int_0^T\tau^\alpha \frac{d}{d\tau}(1+\dm\n w^n(\tau)\dm_2^2)^{-1}\,d\tau\bigg]\,.
\ea
\ee
We estimate the last integral on the right-hand side of \eqref{Rev-Hold}. Integration by parts yields
\[
\ba{l}\displ
-\int_0^T\tau^\alpha \frac{d}{d\tau}(1+\dm\n w^n(\tau)\dm_2^2)^{-1}\,d\tau=-{T^\alpha} (1+\dm\n w^n(T)\dm_2^2)^{-1}\VS\hskip1.7cm+\alpha\int_0^T\tau^{\alpha-1}(1+\dm\n w^n(\tau)\dm_2^2)^{-1}\,d\tau\leq -{T^\alpha} (1+\dm\n w^n(T)\dm_2^2)^{-1}+T^\alpha\leq T^\alpha\,.
\ea
\]
Moreover, employing H\"{o}lder's inequality and recalling the assumptions on $\alpha$ and $\vep$, we have
\[
 \int_0^T\tau^{-\frac{\alpha}{2+\vep}}(1+\dm \n w^n(\tau)\dm_2^2)^{\frac{2}{2+\vep}}\,d\tau\leq \bigg(\int_0^T\tau^{-\frac{\alpha}{\vep}}\,d\tau\bigg)^{\frac{\vep}{2+\vep}}\bigg(\int_0^T(1+\dm \n w^n(\tau)\dm_2^2)\,d\tau\bigg)^{\frac{2}{2+\vep}}\leq c(T,\dm\mathbb P_2w_0\dm_2)\,.
\]
Therefore, going back to \eqref{Rev-Hold}, we conclude the proof of the lemma.
    \ep
    \begin{rem}
        {\rm In the context of the Navier-Stokes equations, in \cite{Mar25} it was proved for the Leray approximating solutions $\{u^n\}_{n\in\N}$ an analogous property to that proved in Lemma\,\ref{le: SECDER-est}, namely
        \[
        \int_0^T \dm D^2 u^n(\tau)\dm_2^\frac23\,d\tau\leq C(T,\dm u_0\dm_2)\,.
        \]
        This property allowed for a strong convergence property related the sequence $\{\n u^n\}_n$, namely
        \be \label{ESt-Dersec-NS}
        \{\n u^n\}_n\to \n u\,\,\text{strongly in }L^p(0,T;L^2(\OO))\,,\,\,\text{for all }p\in [1,2)\,,
        \ee
        where $u$ is the Leray weak solution\footnote{Notice that we have to specify the nature of the weak solution, since estimate \eqref{ESt-Dersec-NS} has not been proved yet for the Hopf-Galerkin approximating sequences in unbounded domains.}. \par In our case, as the norms $\|v(t)\|_\infty$ and $\|\n v(t)\|_3$ have a singularity at $t=0$, we proved the slightly weaker integrability property \eqref{est-DERSEC}.
        However, this property will be sufficient to prove the structure theorem related to the weak solution $w$ to problem \eqref{NS-pert-intro}. In fact, as we will show in the next lemma, it is possible to get a suitable strong convergence property for the sequence $\{\n w^n\}_n$}
    \end{rem}
    \begin{lemma}\label{le: SC-grad}
        {\sl Let $\{w^n\}_n$ be the sequence established in Theorem\,\ref{thm: WP-moll} and let $w$ be the weak solution obtained in Theorem\,\ref{thm: WSOL-per-exist}. Under the assumptions of Lemma\,\ref{le: SECDER-est}, we get
        \[
        \{ \n w^n\}_n\to  \n w\,\,\text{strongly in }L^1(0,T;L^2(\OO))\,.
        \]
        In particular, 
        \[
        \dm \n w^n(t)\dm_2\to \dm\n w(t)\dm_2\,,\,\,\text{a. e. in }(0,T)\,.
        \]}
    \end{lemma}
    \bp
For all $n,m\in \N$, we have
\[
\dm\n (w^n(\tau)-w^m(\tau))\dm_2^2=-(\mathbb P_2\Delta(w^n(\tau)-w^m(\tau)),w^n(\tau)-w^m(\tau))\,.
\]
Employing H\"{o}lder's inequality, we get
\be\label{Hold-SCGrad}
\ba{l}\displ
\int_0^T\dm \n (w^n(\tau)-w^m(\tau))\dm_2\,d\tau\leq \int_0^T \dm \mathbb P_2\Delta(w^n(\tau)-w^m(\tau))\dm_2^\frac12\dm w^n(\tau)-w^m(\tau)\dm_2^\frac12\,d\tau\VS\hskip2.5cm \leq \bigg(\int_0^T \dm \mathbb P_2\Delta(w^n(\tau)-w^m(\tau))\dm_2^\frac{2}{3+\vep}\,d\tau\bigg)^\frac{3+\vep}{4}\bigg(\int_0^T \dm w^n(\tau)-w^m(\tau)\dm_2^{\frac{2}{1-\vep}}\bigg)^\frac{1-\vep}{4}\,.
\ea
\ee
Since
\[
\ba{l}\displ
\dm w^n(\tau)-w^m(\tau)\dm_2\leq \dm w^n(\tau)-w^m(\tau)\dm_{L^2(\OO_R)}+ \dm w^n(\tau)-w^m(\tau)\dm_{L^2(\OO^R)}\VS\hskip3.1cm\leq \dm w^n(\tau)-w^m(\tau)\dm_{L^2(\OO_R)}+ \dm w^n(\tau)\dm_{L^2(\OO^R)}+\dm w^m(\tau)\dm_{L^2(\OO^R)}\,.
\ea
\]
By virtue of Friedrich's Lemma, Lemma\,\ref{le: Conv-extball} and Lemma\,\ref{le: SECDER-est}, we conclude that $\{\n w^n\}_n$ is a Cauchy sequence in $L^1(0,T;L^2(\OO))$. The proof of the lemma can be concluded by classical arguments.
    \ep
    We are now in a position to prove our partial regularity result for the weak solution $w$.
    \begin{tho}\label{Struct-w}
        {\sl Let $w$ be the solution established in Theorem\,\ref{thm: WSOL-per-exist}. Then, for all $\eta>0$, there exist $\theta:= \eta^{-2}(1-\widetilde c\dm\mathbb P_3 v_0\dm_3)^{-1}\dm\mathbb P_2w_0\dm_2^4+1$, and a sequence of open time intervals $\{(\theta_l,T_l)\}_{l\in\N}$ such that $(0,\theta)\setminus \underset{l\in\N}\cup(\theta_l,T_l)$ has zero Lebesgue measure and 
        \begin{equation} \label{reg-w-structhm}
            \ba{c}
            w\in C(\theta_l,T_l;J^{1,2}(\OO))\cap L^2(\theta_l,T_l;W^{2,2}(\OO)),\,\,w_t\in L^2(\theta_l,T_l;L^2(\OO)),\,\forall l\in\N,\VS w\in C([\theta, \infty);J^{1,2}(\OO))\cap L^2(\theta, \infty;W^{2,2}(\OO)),\,\,w_t\in L^2(\theta, \infty;L^2(\OO)). 
            \ea
        \end{equation}}
    \end{tho}
    \bp
By virtue of Lemma\,\ref{le: ENES}, the sequence $\{w^n\}_n$ satisfies the energy inequality
\begin{equation}\label{EN-appr-ST}
    \dm w^n(t)\dm_2^2+2(1-\widetilde c\dm\mathbb P_3 v_0\dm_3)\int_0^t\dm \n w^n(\tau)\dm_2^2\,d\tau\leq \dm \mathbb P_2 w_0\dm_2^2\,.
\end{equation}
We prove that for all $\eta>0$ and for all $n\in\N$ there exists $t^n\in[1,\eta^{-2}(1-\widetilde c\dm\mathbb P_3 v_0\dm_3)^{-1}\dm\mathbb P_2w_0\dm_2^4+1)\equiv[1,\theta)$ such that
\be \label{small-tn} \hskip-0.4cm
\dm \n w^n(t^n)\dm_2\dm\mathbb P_2w_0\dm_2\leq\eta\,.
\ee
If this is not the case, for all $t\in [1,\theta)$, we would have $\eta<\dm \n w^n(t)\dm_2\dm\mathbb P_2w_0\dm_2$. Then
\[
\ba{l}\displ
\dm\mathbb P_2w_0\dm_2=\dm\mathbb P_2w_0\dm_2^4\dm\mathbb P_2w_0\dm_2^{-2}\eta^2\eta^{-2}(1-\widetilde c\dm\mathbb P_3 v_0\dm_3)^{-1}(1-\widetilde c\dm\mathbb P_3 v_0\dm_3)\VS\hskip1.5cm <2\dm\mathbb P_2w_0\dm_2^4\dm\mathbb P_2w_0\dm_2^{-2}\eta^2\eta^{-2}(1-\widetilde c\dm\mathbb P_3 v_0\dm_3)^{-1}(1-\widetilde c\dm\mathbb P_3 v_0\dm_3)\VS\hskip1.6cm=2(1-\widetilde c\dm\mathbb P_3 v_0\dm_3)\eta^2\dm\mathbb P_2 w_0\dm_2^{-2}\int_1^{\eta^{-2}(1-\widetilde c\dm\mathbb P_3 v_0\dm_3)^{-1}\dm\mathbb P_2w_0\dm_2^4+1}dt\VS\hskip2cm<2(1-\widetilde c\dm\mathbb P_3 v_0\dm_3)\int_1^{\eta^{-2}(1-\widetilde c\dm\mathbb P_3 v_0\dm_3)^{-1}\dm\mathbb P_2w_0\dm_2^4+1}\dm\n w^n(t)\dm_2^2 dt\,,
\ea
\]
that contradicts estimate \eqref{EN-appr-ST}. \par We recall that, by virtue of estimates \eqref{IDLrL3-aux} and Corollary\,\ref{cor: LTB-aux}, we have 
\[
\dm v(\tau)\dm_\infty,\dm\n v(\tau)\dm_3\leq c\dm v_0\dm_3\tau^{-\frac{3}{2p}}\,,\,\,\text{for all }\tau>0\,.
\]
Setting $a(\tau):=c(\dm v(\tau)\dm_\infty^2+\dm \n v(\tau)\dm_3^2)$, $\tau>0$, we have $\displ \int_{t^n}^\infty a(\tau)\,d\tau<\infty$.
 By virtue of inequality \eqref{ineq-grad}, setting $y_n(\tau):=\dm \n w^n(\tau)\dm_2^2$, we get
\[
\dot y_n(\tau)\leq a(\tau)y_n(\tau)+cy^3_n(\tau)\,.
\]
Since $y_n(t^n)\leq \frac{\eta^2}{\dm\mathbb P_2w_0\dm_2^2}$ and $\eta$ is allowed to be arbitrarily small, employing Lemma\,\ref{le: Mar-AB} we obtain 
\[
y_n(t)\leq  c(t-t^n)^{-1}\,,\,\,\text{for all }t>t^n\,,
\]
with $c$ independent of $y_n$. Therefore, recalling the relation \eqref{ineq-grad}, we deduce
\[
\int_\theta^\infty \big[\dm D^2 w^n(\tau)\dm_2^2+\dm w^n_\tau(\tau)\dm_2^2\big]\,d\tau\leq C\,,\,\,\text{uniformly in }n\in\N\,.
\]
We conclude that each element of the sequence $\{w^n\}_n$ satisfies
\begin{equation}\label{Regtheta}
    \ba{c}
    w^n\in C(\theta,\infty;J^{1,2}(\OO))\cap L^2(\theta,\infty;W^{2,2}(\OO))\,,\,\,w^n_t\in L^2(\theta,\infty;L^2(\OO))\,,    \ea
\end{equation}
Moreover the bounds detected by the metrics in \eqref{Regtheta} are uniform with respect to $n\in\N$. Hence, we deduce the existence of a subsequence of $\{w^n\}_n$ weakly convergent to a function $\overline w$ in $L^2(\theta,\infty;W^{2,2}(\OO))$, with $\overline w_t\in L^2(\theta,\infty;L^2(\OO))$. By virtue of the uniqueness of the weak limit in $L^2(\theta,\infty;W^{1,2}(\OO))$, we get $\overline w\equiv w$, and, consequently, the weak solution is regular in $[\theta,\infty)$. Therefore, there exists a pressure field $\pi_w$ such that $(w,\pi_w)$ is a regular solution to \eqref{NS-pert-intro} on $[\theta,\infty)\times\OO$. \par We now focus on the regularity of the weak solution on the interval of time $(0,\theta)$. By virtue of Lemma\,\ref{le: SC-grad}, we have
\be\label{CQO}
\dm \n w^n(t)\dm_2\to \dm \n w(t)\dm_2\,,\,\,\text{a. e. in }(0,\theta)\,.
\ee
We set
\[
\mathscr I:=\{t\in (0,\theta)\,:\, \text{the limit property \eqref{CQO} holds}\}\,. 
\]
Let $t_\alpha\in (0,\theta)\cap \mathscr I$. We recall that the following inequality holds:
\be \label{DiffIneq-theta}
\dot y_n(\tau)\leq a(\tau)y_n(\tau)+cy^3_n(\tau)\,,\,\,\text{for all }\tau >t_\alpha\,.
\ee
Moreover, as $t_\alpha\in \mathscr I$, there exists $n(t_\alpha)\in\N$ such that, for all $n\ge n(t_\alpha)$, $y_n(t_\alpha)\leq \dm \n w(t_\alpha)\dm_2+1$. Therefore, for $n\ge n(t_\alpha)$, the inequality \eqref{DiffIneq-theta} can be integrated over a maximal interval of time $[t_\alpha,t_\alpha+T(n,t_\alpha))\supset \bigg[t_\alpha, t_\alpha+\frac{1}{2c}\frac{e^{-\int_{t_\alpha}^\tau a(s)\,ds}}{(\dm \n w(t_\alpha)\dm_2+1)^2}\bigg)$.
As a matter of fact, setting $z_n(\tau):=y_n^{-2}(\tau)$, from inequality \eqref{DiffIneq-theta}, we get
\[
z_n(\tau)\geq e^{-\int_{t_\alpha}^\tau a(s)\,ds}\bigg[z_n(t_\alpha)-e^{\int_{t_\alpha}^\tau a(s)\,ds}\int_{t_\alpha}^\tau 2c\,ds\bigg]\,.
\]
Since $z_n(t_\alpha)\geq  (\dm \n w(t_\alpha)\dm_2+1)^{-2}$ for all $n\geq n(t_\alpha)$, then we easily deduce
\[
[t_\alpha,t_\alpha+T(n,t_\alpha))\supset \bigg[t_\alpha, t_\alpha+\frac{1}{2c}\frac{e^{-\int_{t_\alpha}^\tau a(s)\,ds}}{(\dm \n w(t_\alpha)\dm_2+1)^2}\bigg)\,,\,\,\text{for all }n\geq n(t_\alpha)\,.
\]
Therefore, we have that
\[
[t_\alpha,t_\alpha+T(t_\alpha)):=\underset{n\ge n(t_\alpha)}\cap [t_\alpha,t_\alpha+T(n,t_\alpha))\,\,\text{is not empty}
\]
and we obtain
\be \label{reg-0theta}
\dm \n w^n(s)\dm_2^2+\int_{t_\alpha}^s \big[\dm D^2 w^n(\tau)\dm_2^2+\dm w^n_\tau(\tau)\dm_2^2\big]\,d\tau<\infty\,,\,\,\text{for all }s\in [t_\alpha,t_\alpha+T(t_\alpha))\,.
\ee
  By virtue of the bound \eqref{reg-0theta} we state the existence of a subsequence of $\{w^n\}_n$ that admits a limit $\overline{\overline w}$, regular solution to \eqref{NS-pert-intro} on some maximal interval of time $I_\alpha:=[t_\alpha,t_\alpha+T(t_\alpha))\supset \bigg[t_\alpha, t_\alpha+\frac{1}{2c}\frac{e^{-\int_{t_\alpha}^\tau a(s)\,ds}}{(\dm \n w(t_\alpha)\dm_2+1)^2}\bigg)$. As the weak limit in $L^2(0,\theta;W^{1,2}(\OO))$ is unique, we get $\overline {\overline w}\equiv w$. Therefore, the weak solution $w$ is also regular on $I_\alpha\times\OO$. \par We now consider $t_\beta\in (0,\theta)\cap (\mathscr I\setminus I_\alpha)$. By analogous arguments, we get that the weak solution is regular on some maximal interval of time $I_\beta=[t_\beta,t_\beta+T(t_\beta))$. One of the following holds:
\begin{description}
    \item[i.] $I_\alpha\cap I_\beta=\emptyset $;
    \item[ii.] $I_\alpha\cap I_\beta\ne \emptyset $, then $I_\alpha\subset I_\beta$.
\end{description}
We justify the second item. If $I_\alpha\cap I_\beta\ne\emptyset$, since $t_\beta\notin I_\alpha$, then $t_\beta<t_\alpha$ and, since $I_\beta$ is a maximal interval of existence, we have $t_\alpha\in I_\beta$. Therefore, it has to be $I_\alpha\subset I_\beta$. If item $\bf i.$ holds, we take $t_\gamma\in (0,\theta)\cap (\mathscr I\setminus (I_\alpha\cup I_\beta))$. On the other hand, if item $\bf ii.$ holds, we replace $I_\alpha$ by $I_\beta$ and take $t_\gamma\in (0,\theta)\cap (\mathscr I\setminus I_\beta)$. Iterating this procedure, we construct a family $\{I_\alpha\}_{\alpha\in\mathbb A}$ of time intervals of regularity for the weak solution $w$ such that
\be \label{incl-struc}
(0,\theta)\cap \mathscr I\subseteq \underset{\alpha\in \mathbb A}\cup I_\alpha\,. 
\ee
Since for all $\alpha\in \mathbb A$ the interior set of $I_\alpha$ (denoted by $\mathring I_\alpha$) is not empty and, if $\beta\ne \alpha$, $\mathring I_\alpha\cap \mathring I_\beta=\emptyset$, then $\mathbb A$ is at most countable. Moreover, the inclusion property \eqref{incl-struc} ensures that
\[
(0,\theta)\setminus \underset{\alpha\in\mathbb A}\cup \mathring I_\alpha
\]
has zero Lebesgue measure. The theorem is completely proved.
    \ep
    \bp[Proof of Theorem\,\ref{thm: WSOL-pert}]
The proof of the Theorem follows from Theorem\,\ref{thm: WSOL-per-exist} and Theorem\,\ref{Struct-w}.
    \ep
    \section{Proof of Theorem\,\ref{thm: MR}}\label{sec: PMR}
    We give the proof of our main result.
    \bp[Proof of Theorem\,\ref{thm: MR}]
We divide the proof in two steps: we first prove the existence statement and then the structure theorem.\par 
{\bf Step 1.} (\textit{Existence}). Following Lemma\,\ref{Dec-Lp}, we decompose $u_0$ as the sum
\[
u_0=v_0+w_0\,,
\]
with $v_0\in L^3(\OO)\cap L^p(\OO)$ and $w_0\in L^2(\OO)\cap L^p(\OO)$.
We denote by $c(3)$ the constant related to estimate \eqref{HD-ineq} with $p=3$. Choosing $\rho=\bigg(\frac{\overline\xi_0}{2 c(3)\dm u_0\dm_p^\frac{p}{3}}\bigg)^{\frac{3}{3-p}}$ in Lemma\,\ref{Dec-Lp}, invoking \eqref{HD-ineq}, we have
\[
\dm\mathbb P_3 v_0\dm_3<\overline \xi_0\,.
\]
Therefore, we can construct our solution $u$ as the sum $u=v+w$, where $v$ is solution to problem \eqref{NS-aux-intro} with initial datum $v_0\in L^3(\OO)\cap L^p(\OO)$, whose existence is ensured by Theorem\,\ref{thm: L3-weakdata}, while $w$ is a weak solution to problem \eqref{NS-pert-intro} with initial datum $w_0\in L^2(\OO)\cap L^p(\OO)$, ensured by Theorem\,\ref{thm: WSOL-pert}. We show that $u$ satisfies the items $\bf (a)-(d)$ of Definition\,\ref{def: WS-NS}. By virtue of the regularity properties \eqref{reg-sol} for $v$ and the ones given in Definition\,\ref{def: WS-pert} for $w$, we deduce
\[
u\in L^4(0,T;L^3(\OO))\,.
\]
Moreover, given $\OO'$ a bounded subdomain of $\OO$ by straightforward considerations we deduce
\[
v\in L^{2,\infty}(0,T;W^{1,2}(\OO'))\,.
\]
Therefore, recalling the regularity properties of $w$, we get
\[
u\in L^{2,\infty}(0,T;W^{1,2}_{loc}(\OO))\,.
\]
Items $\bf b)-c)$ and the trace property in $\bf d)$ for $u$ are also easy to verify, since they hold for $w$ and $v$ is a global regular solution satisfying properties \eqref{reg-sol}. Concerning the convergence to the initial datum, we see that
\[
(u(t),\varphi)=(v(t),\varphi)+(w(t),\varphi)\to (v_0,\varphi)+(w_0,\varphi)=(u_0,\varphi)\,,\,\,\text{as }t\to 0\,,\,\text{ for all }\varphi\in\mathscr C_0(\OO)\,.
\]
The existence of a weak solution in the sense of Definition\,\ref{def: WS-NS} is achieved.\par
{\bf Step 2.} (\textit{Regularity}). We have $u=v+w$ and we know $v$ satisfies the regularity properties expressed in Definition\,\ref{def: SOLR-NS} for all times. Concerning $w$, in Theorem\,\ref{Struct-w} we proved the regularity properties \eqref{reg-w-structhm}. Therefore, we easily conclude that $u$ has the regularity properties detected in Definition\,\ref{def: SOLR-NS} in all the time intervals in which the regularity of $w$ is ensured. The theorem is completely proved.  
    \ep
    \bp[Proof of Corollary\,\ref{cor: LTB-sol}]
Since the solution constructed in Theorem\,\ref{thm: MR} has the regularity properties detected in Lemma\,\ref{le: GWP-L3} on the interval of time $(\theta,\infty)$, then properties \eqref{IDLrL3} are satisfied for $t>\theta$. Therefore, recalling also Corollary\,\ref{cor: LTB-aux}, the claim of the corollary follows.
    \ep
\section*{Appendix: On the proof of some properties stated in Lemma\,\ref{le: GWP-L3}}
Let $\OO\subseteq \R^3$ be the whole space, a half-space or an exterior domain with a sufficiently smooth boundary.
As already remarked in the main body of the paper, the proof of the existence of a global regular solution to problem \eqref{NS-aux-intro} for an initial datum in $J^3(\OO)$, satisfying a suitable smallness condition, is known in the literature thanks to works \cite{Giga86,Iwa, Kato, Mar-DCDS}. \par However, in the case of $\OO$ exterior domain, for the sake of completeness we find it worth to spend some words on the proof of properties \eqref{reg-sol-prel}$_{4,5}$ and \eqref{IDLrL3}. In fact, the case of the exterior problem is particularly interesting since the semigroup properties of the Stokes operator are different with respect to the Cauchy problem and the half-space problem (see Lemma\,\ref{le: Stokes-Jp} and Remark\,\ref{rem: Cauchy+semisp-Semi}).\par We first recall a result related to the Stokes initial boundary value problem with initial data in $L^1(\OO)$, proved in \cite{Mar-L1}. Let us consider the following problem:
\begin{equation}\label{Sto-L1}
\begin{cases}
V_t-\Delta V+\n\pi_V=f\,,\,\,\quad\quad\text{in }(0,T)\times\OO,\\
\hskip1.83cm\n\cdot V=0\,,\,\,\quad\quad\text{in }(0,T)\times\OO\,,\\
\hskip1.6cmV(t,x)=0\,,\,\,\,\,\,\,\quad\,\,\text{in }(0,T)\times\partial\OO\,,\\
\hskip0.72cm(V(0,x),\varphi)=(V_0(x),\varphi)\,,\,\,\,\text{for all }\varphi\in\mathscr C_0(\Omega)\,.
\end{cases}
\end{equation}
Let $r\in (1,\infty)$ and let $r'$ be the complementary exponent of $r$. We set
\[
W_{r'}:=\{\phi\,:\, \phi\in C^1([0,T]\times\overline \OO)\cap C([0,T];J^{1,r'}(\OO))\,,\,\,\phi_t\in C([0,T];L^{r'}(\OO))\};
\]
\begin{lemma}\label{le: Sto-L1}
{\sl Let $V_0\in L^1(\OO)$. There exists a unique solution $(V,\pi_V)$ to problem \eqref{Sto-L1} such that
\begin{description}
\item[i.] For all $\eta>0$, $r>1$,
$$
V\in C(\eta,T;J^r(\OO))\cap L^\infty(\eta,T;J^{1,r}(\OO)\cap W^{2,r}(\OO)),\,\n\pi_V\in L^\infty(\eta,T;L^r(\OO));
$$
    \item[ii.] \be \label{sem-L1}
    \ba{rl}\displ
    \dm V(t)\dm_r&\hskip-0.2cm\leq c\dm V_0\dm_1t^{-\mu}\,,\quad \mu=\frac32\big(1-\frac1r\big)\,,\,\,t>0\,,\,r\in (1,\infty];\VS \dm \n V(t)\dm_r&\hskip-0.2cm\leq c\dm V_0\dm_1t^{-\mu_1}\,,\,\,\,\, \mu_1=\begin{cases}
        \frac12+\mu \,,\,\,\text{if }t\in(0,1)\,,\,r\in (1,\infty)\,,\\ \frac12+\mu\,,\,\,\text{if }t>0\,\,\text{and }r\in (1,3]\,,\\ \frac{3}{2}\,,\quad\,\,\,\,\,\text{if }t\ge1\,\,\text{and }r> 3\,;
    \end{cases}\VS
    \dm V_t(t)\dm_r&\hskip-0.2cm\leq c\dm V_0\dm_1t^{-\mu_2}\,,\,\,\,\,\mu_2=1+\mu\,,\,\,t>0\,,\,r\in (1,\infty]\,;
    \ea
    \ee 
    \item[iii.] $\displ \int_0^t \big[(V(\tau),\phi(\tau))-(\n V(\tau),\n\phi(\tau))\big]\,d\tau =(V(t),\phi(t))-(V_0,\phi(0))$, for all $\phi\in W_{q'}$, provided that $\frac12+\frac32(1-\frac1q)<1$;
    \item[iv.] $\displ \lim_{t\to 0}(V(t),\varphi)=(V_0,\varphi)$, for all $\varphi\in\mathscr C_0(\OO)$.
    \end{description}}
\end{lemma}
 In order to prove properties \eqref{reg-sol-prel}$_4$ and \eqref{IDLrL3}$_2$ we are going to employ the semigroup formalism, as done by Kato in \cite{Kato}, while to prove \eqref{IDLrL3}$_1$ we employ a duality argument. In order to employ the semigroup formalism, we recall that for $v_0\in J^3(\OO)$ problem \eqref{NS-aux-intro} can be rewritten in the abstract integral formulation 
\be \label{AIP}
v(t)=e^{-tA}v_0+\int_0^t e^{-(t-\tau)A}\mathbb P_3[v(\tau)\cdot\n v(\tau)]\,d\tau\,,
\ee
where we denoted by $A$ the Stokes operator on $J^3(\OO)$ and we considered the kinetic field $v$ as a function of $t$ taking values in a Banach space. \par We recall that the solution $v$ is obtained by the recursive approximation technique. Hence, it is deduced as the limit, in suitable metrics, of the sequence $\{v^n\}_{n\in\N_0}$, with
\be \label{rec-appr}
\ba{c}\displ
v^0(t):=e^{-tA}v_0\,, \,\,v^{n+1}(t):=v^0(t)+\int_0^te^{-(t-\tau)A}\mathbb P_3[v^n(\tau)\cdot\n v^n(\tau)]\,d\tau\,.
\ea
\ee
\bp[Proof of properties \eqref{reg-sol-prel}$_{4,5}$]
By virtue of Lemma\,\ref{le: Stokes-Jp} we have
\be \label{iter-Sto}
t^\frac12\n v^0(t)\in C([0,\infty);L^3(\OO))\,.
\ee
Let $n\in \N$ and assume that
\[
t^\frac12\n v^n(t)\in C([0,\infty);L^3(\OO))\,,\,\,\text{with a bound independent of }n\,.
\]
We prove by induction that 
\[
t^\frac12\n v^{n+1}(t)\in C([0,\infty);L^3(\OO))\,,\,\,\text{with a bound independent of }n\,.
\]
By virtue of \eqref{rec-appr} we have
\[
\dm \n v^{n+1}(t)\dm_3\leq \dm \n v^0(t)\dm_3+\int_0^t \dm \n e^{-(t-\tau)A}\mathbb P_3[v^n(\tau)\cdot\n v^n(\tau)]\dm_3\,d\tau.
\]
As \eqref{iter-Sto} holds, we only have to discuss the integral term. Let $\delta\in (0,1)$. Employing estimate \eqref{Sem-prop}$_2$ and recalling \eqref{HD-ineq}, we obtain
\[
\int_0^t \dm \n e^{-(t-\tau)A}\mathbb P_3[v^n(\tau)\cdot\n v^n(\tau)]\dm_3\,d\tau\leq c\int_0^t (t-\tau)^{-\frac12 - \frac32 (\frac{1+\delta}{3}-\frac13)}\dm v^n(\tau)\cdot\n v^n(\tau)\dm_{\frac{3}{1+\delta}}\,d\tau\,.
\]
Employing H\"{o}lder's inequality we get
\[
\int_0^t (t-\tau)^{-\frac12 - \frac32 (\frac{1+\delta}{3}-\frac13)}\dm v^n(\tau)\cdot\n v^n(\tau)\dm_{\frac{3}{1+\delta}}\,d\tau\leq \int_0^t (t-\tau)^{-\frac12 - \frac32 (\frac{1+\delta}{3}-\frac13)}\dm v^n(\tau)\dm_{\frac3\delta}\dm\n v^n(\tau)\dm_{3}\,d\tau\,.
\]
Finally, by virtue of property \eqref{ASYMP-SOLGLO-AUX}$_1$ and the induction hypothesis, we conclude that
\[
\dm \n v^{n+1}(t)\dm_3\leq \dm \n v^0(t)\dm_3+c\dm v_0\dm_3^2\int_0^t(t-\tau)^{-\frac12 - \frac32 (\frac{1+\delta}{3}-\frac13)}\tau^{-\frac12-\frac32(\frac13-\frac\delta3)}\,d\tau\,,
\]
and
\[
\int_0^t(t-\tau)^{-\frac12 - \frac32 (\frac{1+\delta}{3}-\frac13)}\tau^{-\frac12-\frac32(\frac13-\frac\delta3)}\,d\tau=t^{-\frac12 }\int_0^1 s^{-1+\frac\delta2 }(1-s)^{-\frac12-\frac\delta2}\,ds=ct^{-\frac12}\,.
\]
We have thus proved
\[
t^\frac12\n v^{n+1}(t)\in C([0,\infty);L^3(\OO))\,.
\]
Therefore, for all $n\in \N$,
\[
t^\frac12\n v^{n}(t)\in C([0,\infty);L^3(\OO))\,.
\]
We point out that the smallness condition on $\dm v_0\dm_3$ allows us to get the aforementioned property uniformly with respect to $n\in\N$ and this allows us to state 
\[
t^\frac12 \n v\in C([0,\infty);L^3(\OO))\,.
\]
Concerning property \eqref{reg-sol-prel}$_5$ it can be proved by following analogous arguments. However, we have to restrict to the interval of time $[0,1)$ in view of the semigroup property \eqref{Sem-prop}$_2$.
\ep
We recall that \eqref{rec-appr}$_2$ is the integral formulation of the following  initial boundary value problem
\begin{equation}\label{PDE-iter}
    \begin{cases}
        v^{n+1}_t-\Delta v^{n+1}+v^n\cdot\n v^n+\n \pi^{n+1}=0\,,\,\,\text{in }(0,T)\times\OO\,,\\ \hskip2cm\n\cdot v^{n+1}=0\,,\qquad\qquad\qquad\,\,\text{in }(0,T)\times\OO\,,\\ \hskip2.6cm v^{n+1}=0\,,\qquad\qquad\qquad\,\,\text{in }(0,T)\times\partial\OO\,,\\ \hskip1.73cm v^{n+1}(0,x)=0\,,\qquad\qquad\qquad\,\,\text{in }\OO\,.
    \end{cases}
\end{equation}
\bp[Proof of property \eqref{IDLrL3}$_1$]
Employing a duality argument, we are going to prove that
\be \label{v-3intersr}
\dm v(t)\dm_r \leq c\dm v_0\dm_r\,,\,\, \dm v(t)\dm_\infty\leq c\dm v_0\dm_rt^{-\frac{3}{2r}}\,,
\ee
so that, employing the Riesz-Thorin theorem, we arrive at \eqref{IDLrL3}$_1$. \par We proceed by induction on the recursive approximations. It is readily seen that
\[
\dm v^0(t)\dm_r\leq c\dm v_0(t)\dm_r\,,\,\,\text{for all }t>0\,.
\]
Let $n\in \N$ and assume that 
\[
\dm v^n(t)\dm_r\leq c\dm v_0(t)\dm_r\,,\,\,\text{for all }t>0\,,\,\,\text{with c independent of }n\,.
\]
We prove
\be \label{rec-rinter3}
\dm v^{n+1}(t)\dm_r\leq c\dm v_0(t)\dm_r\,,\,\,\text{for all }t>0\,,\,\,\text{with c independent of }n\,.
\ee
We consider $(z,\pi_z)$ solution to the Stokes problem with initial datum $z_0\in \mathscr C_0(\Omega)$. Let $t>0$. We set $\widehat z(\tau ,x)= z(t-\tau,x)$, $\tau \in (0,t)$. We multiply equation \eqref{PDE-iter}$_1$ by $\widehat z$ and integrate over $(0,t)\times \Omega$. We get
\be \label{dual-uniq}
(v^{n+1}(t),z_0)=\int_0^t ( v^{n}\cdot \n \widehat z, v^n)\,d\tau \,.
\ee
Using H\"{o}lder's inequality, the semigroup properties for $z$, property \eqref{ASYMP-SOLGLO}$_1$ and the induction hypothesis, we find that
\[
\ba{l}\displ
\Bigg|\int_0^t ( v^{n}\cdot \n \widehat z, v^{n})\,d\tau\Bigg|\leq \int_0^t \dm v^n(\tau)\dm_\infty\dm \n \widehat z(\tau)\dm_{\frac{r}{r-1}}\dm v^n(\tau)\dm_r\,d\tau\VS \hskip4cm\leq c\dm v_0\dm_r\dm v_0\dm_3\dm z_0\dm_{\frac{r}{r-1}}\int_0^t(t-\tau)^{-\frac12}\tau^{-\frac12}\,d\tau\,,
\ea
\]
with the last integral being clearly finite. The smallness condition on $\dm v_0\dm_3$ and the arbitrariness of $z_0$ allow us to deduce the validity of \eqref{rec-rinter3}. Hence, we get property \eqref{v-3intersr}$_1$. \par To prove \eqref{v-3intersr}$_2$, by using analogous arguments and also employing \eqref{sem-L1}$_2$, we point out that
\[
\ba{l}\displ
\Bigg|\int_0^t ( v^{n}\cdot \n \widehat z, v^{n})\,d\tau\Bigg|\leq \int_0^t \dm v^n(\tau)\dm_\infty\dm \n \widehat z(\tau)\dm_{\frac{2r}{2r-1}}\dm v^n(\tau)\dm_{2r}\,d\tau\VS \hskip4cm\leq c\dm v_0\dm_r\dm v_0\dm_3\dm z_0\dm_{1}\int_0^t(t-\tau)^{-\frac12-\frac{3}{4r}}\tau^{-\frac12-\frac{3}{4r}}\,d\tau\,.
\ea
\]
Since $r>\frac32$, we have
\[
\int_0^t(t-\tau)^{-\frac12-\frac{3}{4r}}\tau^{-\frac12-\frac{3}{4r}}\,d\tau=t^{-\frac{3}{2r}}\int_0^1 (1-s)^{-\frac12-\frac{3}{4r}}s^{-\frac12-\frac{3}{4r}}\,ds=ct^{-\frac{3}{2r}}\,.
\]
Hence, in view of the smallness condition on $\dm v_0\dm_3$, we obtain, uniformly in $n\in\N$,
\[
\dm v^n(t)\dm_\infty\leq c\dm v_0\dm_rt^{-\frac{3}{2r}}\,.
\]
Therefore, \eqref{v-3intersr}$_2$ follows. Employing the Riesz-Thorin theorem for the identity operator, we get \eqref{IDLrL3}$_1$. 
\ep
\bp[Proof of property \eqref{IDLrL3}$_2$ (Exterior domain case)] 
We employ again the abstract integral formulation and the properties of the Stokes semigroup. We prove \eqref{IDLrL3}$_2$ first for $q=r$, then for $q=3$, and employ the Riesz-Thorin theorem to conclude. \par Let us prove \eqref{IDLrL3}$_2$ for $q=r$. We employ an induction argument on the recurrence approximations. It is readily seen that
\[
\dm \n v^0(t)\dm_r\leq c\dm v_0\dm_rt^{-\frac12}\,,\,\,\text{for all }t>0\,.
\]
Let $n\in\N$ and assume that
\[
\dm \n v^n(t)\dm_r\leq c\dm v_0\dm_rt^{-\frac12}\,,\,\,\text{for all }t>0\,\,\text{uniformly in }n\in\N\,.
\]
We prove
\be \label{iter-gradLr}
\dm \n v^{n+1}(t)\dm_r\leq c\dm v_0\dm_rt^{-\frac12}\,,\,\,\text{for all }t>0\,\,\text{uniformly in }n\in\N\,.
\ee
We have
\[
\dm \n v^{n+1}(t)\dm_r\leq \dm \n v^0(t)\dm_r+\int_0^t \dm \n e^{-(t-\tau)A}\mathbb P_3[v^n(\tau)\cdot\n v^n(\tau)]\dm_r\,d\tau\,. 
\]
Let $\delta\in (0,\frac{r}{3})$. Employing estimates \eqref{Sem-prop}$_2$ and \eqref{HD-ineq}, we find that
\[
\int_0^t \dm \n e^{-(t-\tau)A}\mathbb P_3[v^n(\tau)\cdot\n v^n(\tau)]\dm_r\,d\tau\leq \int_0^t (t-\tau)^{-\frac12-\frac{3\delta}{2r}}\dm v^n(\tau)\cdot\n v^n(\tau)\dm_{\frac{r}{1+\delta}}\,d\tau\,.
\]
By virtue of H\"{o}lder's inequality, property \eqref{ASYMP-SOLGLO}$_1$ and the induction hypothesis, we get
\[
\ba{l}\displ
\int_0^t (t-\tau)^{-\frac12-\frac{3\delta}{2r}}\dm v^n(\tau)\cdot\n v^n(\tau)\dm_{\frac{r}{1+\delta}}\,d\tau\leq \int_0^t (t-\tau)^{-\frac12-\frac{3\delta}{2r}}\dm v^n(\tau)\dm_{\frac{r}{\delta}}\dm\n v^n(\tau)\dm_r\,d\tau\VS \hskip0.5cm \leq c\dm v_0\dm_3\dm v_0\dm_r\int_0^t (t-\tau)^{-\frac12-\frac{3\delta}{2r}}\tau^{-1+\frac{3\delta}{2r}}\,d\tau=c\dm v_0\dm_3\dm v_0\dm_rt^{-\frac12}\int_0^1 s^{-1+\frac{3\delta}{2r}}(1-s)^{-\frac12-\frac{3\delta}{2r}}\,ds\,,
\ea
\]
the last integral being finite. Invoking the smallness condition on $\dm v_0\dm_3$, we deduce \eqref{iter-gradLr}. Therefore, for all $n\in\N$, we conclude that
\[
\dm \n v^n(t)\dm_r\leq c\dm v_0\dm_rt^{-\frac12}\,,\,\,\text{for all }t>0\,,\,\,\text{c independent of }n\in\N\,,
\]
and, therefore,
\[
\dm \n v(t)\dm_r\leq c\dm v_0\dm_rt^{-\frac12}\,,\,\,\text{for all }t>0\,.
\]
By analogous arguments, we are able to prove that
\[
\dm \n v(t)\dm_3\leq c\dm v_0\dm_rt^{-\frac32(\frac1r-\frac13)-\frac12}\,,
\]
and, employing the Riesz-Thorin theorem for the identity operator, we conclude the proof of \eqref{IDLrL3}$_2$. 
\ep

{\bf Acknowledgments. }{\small I would like to thank Prof. P. Maremonti for his valuable supervision of this research. I am also thankful to Prof. Z. Bradshaw for his comments and for bringing to my attention some further references.}

{\bf Data availability statement. }{\small Data sharing not applicable to this article as no datasets were
generated or analyzed during the current study.}

{\bf Conflict of interest. }{\small The author states
that there is no Conflict of interest.}

\end{document}